\documentclass[lfeqn,12pt]{article}
\usepackage{amssymb,bbm,pifont,mathtools,mathabx,mathrsfs,%fdsymbol,
}
\usepackage{minitoc}
 
\usepackage[nottoc]{tocbibind}

\usepackage{comment}
\usepackage{multirow}
\usepackage[multiple]{footmisc}
\usepackage{xcolor}
\usepackage[colorlinks=true]{hyperref}
\usepackage{marginnote}
\newcommand{\Giu}{{\bigskip\noindent}}

\newcommand{\noi}{{\noindent}}

\newtheorem{theorem}{Theorem}
\newtheorem{definition}[theorem]{Definition}
\newtheorem{proposition}[theorem]{Proposition}
\newtheorem{lemma}[theorem]{Lemma}
\newtheorem{remark}[theorem]{Remark}
\newtheorem{sublemma}[theorem]{Sublemma}
\newtheorem{corollary}[theorem]{Corollary}
\newtheorem{assumption}[theorem]{Assumption}
\newtheorem{question}[theorem]{Question}
\newtheorem{conjecture}[theorem]{Conjecture}
\newtheorem{notationalrem}[theorem]{Notational Remark}
\newtheorem{tools}[subsection]{$\negsp\negsp$}
\newcommand\asm[1]{ \begin{assumption}\label{#1} }
\newcommand\easm{ \end{assumption} }
\newcommand\qst[1]{ \begin{question}\label{#1} }
\newcommand\eqst{ \end{question} }
\newcommand\conj[1]{ \begin{conjecture}\label{#1} }
\newcommand\econj{ \end{conjecture} }
\newcommand\dfn[1]{ \begin{definition}\label{#1} }
\newcommand\dfntwo[2]{ \begin{definition}[#2]\label{#1} }
\newcommand\edfn{ \end{definition} }
\newcommand\rem[1]{ \begin{remark}\label{#1} \small \rm}
\newcommand\remtwo[2]{ \begin{remark}[#2]\label{#1} \rm}
\newcommand\erem{ \end{remark} }
\newcommand\thm[1]{ \begin{theorem}\label{#1}}
\newcommand\thmtwo[2]{ \begin{theorem}[#2]\label{#1}}
\newcommand\ethm{ \end{theorem} }
\newcommand\pro[1]{ \begin{proposition}\label{#1}}       
\newcommand\protwo[2]{ \begin{proposition}[#2]\label{#1}}
\newcommand\epro{ \end{proposition} }
\newcommand\lem[1]{ \begin{lemma}\label{#1}}
\newcommand\lemtwo[2]{ \begin{lemma}[#2]\label{#1}}
\newcommand\elem{ \end{lemma} }
\newcommand\sublem[1]{ \begin{sublemma}\label{#1}}
\newcommand\sublemtwo[2]{ \begin{sublemma}[#2]\label{#1}}
\newcommand\esublem{ \end{sublemma} }
\newcommand\cor[1]{ \begin{corollary}\label{#1}}
\newcommand\cortwo[2]{ \begin{corollary}[#2]\label{#1}}
\newcommand\ecor{ \end{corollary} }
\newcommand\notrem[1]{ \begin{notationalrem}\label{#1} \sl}
\newcommand\enotrem{ \end{notationalrem} }

\newcommand\average[1]{{ \left\langle #1 \right\rangle}}

\newcommand\equ[1]{{\rm (\ref{#1})}}
\newcommand\beq[1]{ \begin{equation}\label{#1} }
\newcommand{\eeq}{ \end{equation} }

\newcommand\beqa[1]{ \begin{eqnarray} \label{#1}}
\newcommand{\eeqa}{ \end{eqnarray} }
\newcommand{\beqano}{ \begin{eqnarray*} }
\newcommand{\eeqano}{ \end{eqnarray*} }

\newcommand{\proof}{\par\medskip\noindent{\bf Proof\ }}

\newcommand{\ie}{{\it i.e.\  }}
\newcommand{\eg}{{\it e.g. \ }}

\newcommand{\ifof}{{\it iff\ }}

\newcommand{\dst}{\displaystyle}
\newcommand{\qed}{\hskip.5truecm
            \vrule width 1.7truemm height 3.5truemm depth 0.truemm
            \par\Giu}

\newcommand\ovl[1]{ \overline {#1} }

\newcommand\su[1]{ \frac{1}{ {#1}} }
\newcommand\minfoc{ {\, \rm minfoc\, }}

\newcommand\supp{ {\, \rm supp\, }}

\newcommand\meas{ {\, \rm meas\, }}

\newcommand{\diag}{{ \, \rm diag \, }}

\newcommand{\Id}{ {\rm Id }}

\newcommand{\dpr}{ {\partial}   }

\newcommand\eqby[1]{\stackrel{\equ{#1}}{=}}
\newcommand\leby[1]{\stackrel{\equ{#1}}{\le}}
\newcommand\ltby[1]{\stackrel{\equ{#1}}{<}}

\renewcommand{\Im}{{\rm \, Im\,}}
\renewcommand{\Re}{{\rm \, Re\,}}

\newcommand{\negsp}{\hspace{-.04truecm}}
\newcommand\ex{\, e}%{\mathbbmss e}}

\renewcommand{\a }{ {\alpha}   }
\renewcommand{\b}{ {\beta}   }
\newcommand{\g}{ {\gamma}   }
\newcommand{\G}{ {\Gamma}   }
\renewcommand{\d}{ {\delta}   }
\newcommand{\D}{ {\Delta}   }
\newcommand{\vae }{ {\varepsilon}   }
\newcommand{\z}{ {\zeta} }
\renewcommand{\eta}{ {\theta}   }

\renewcommand{\k}{ {\kappa}   }
\renewcommand{\l}{ {\lambda}   }

\newcommand{\n}{ {\nu}   }
\newcommand{\x }{ {\xi}   }

\newcommand{\p}{ {\pi}   }

\renewcommand{\r}{ {\rho}   }
\newcommand{\s}{ {\sigma}   }

\renewcommand{\t}{ {\tau}   }
\newcommand{\f}{ {\varphi}   }

\renewcommand{\o}{ {\omega}   }

\newcommand{\torus}{ {\mathbb{ T}}   }
\renewcommand{\natural}{ {\mathbb{ N}}   }
\newcommand{\real}{ {\mathbb{ R}}   }
\newcommand{\integer}{ {\mathbb{ Z}}   }
\newcommand{\complex}{ {\mathbb { C}}   }

\newcommand{\tn}{ {\torus^d} }

\newcommand{\rn}{ {\real^d}   }

\newcommand{\cn}{ {\complex^d }   }
\newcommand{\zn}{ {\integer^d }   }
\newcommand{\nn}{ {\natural^d }   }

\newcommand\ppu{{ (1) }}
\newcommand\ppd{{ (2) }}
\newcommand\ppt{{ (3) }}

\font\teneufm=eufm10
\font\seveneufm=eufm7
\font\fiveeufm=eufm5
\newfam\eufmfam
\textfont\eufmfam=\teneufm
\scriptfont\eufmfam=\seveneufm
\scriptscriptfont\eufmfam=\fiveeufm

\newcommand\appA[1]{\section{#1}\label{app:A}
\renewcommand{\theequation}{A.\arabic{equation}}
           \setcounter{equation}{0}
\renewcommand{\thetheorem}{A.\arabic{theorem}}
           \setcounter{theorem}{0}
                  }
\newcommand\appB[1]{\section{#1}%\label{app:B}
\renewcommand{\theequation}{B.\arabic{equation}}
           \setcounter{equation}{0}
\renewcommand{\thetheorem}{B.\arabic{theorem}}
           \setcounter{theorem}{0}           
           }
\newcommand\appC[1]{\section{#1}%\label{app:E}
\renewcommand{\theequation}{C.\arabic{equation}}
           \setcounter{equation}{0}
\renewcommand{\thetheorem}{C.\arabic{theorem}}
           \setcounter{theorem}{0}
           }
\def\uno{{\mathbbm 1}}

\def\id{{\rm id }}
\newcommand{\wh}{\widehat}
\newcommand{\wt}{\widetilde}

\begin{document}

%%%%%%%%%%%%%%%%%%%%%%%%%%%
%%%%%%%   title page

\date{\small \today}

\title{{\bf  
A KAM Theorem for finitely differentiable Hamiltonian systems%\\Persistence of a Cantor--like familly of KAM tori under a priori optimal regularity assumption\\
%Integrability of Hamiltonians on Cantor--like sets under  a priori sharp regularity assumption
}%%\thanks{{\bf Acknowledgments.} We are grateful to Prof. Luigi Chierchia for useful discussions. The work of this paper was mainly done while the author was visiting the Abdus Salam International Centre for Theoretical Physics in Trieste, Italy. 
%%%We would like to thank Prof. L. Chierchia and Prof. L. Biasco for helping a lot and their patience.
%%%work was completed while both authors were affiliated to the Abdus Salam International Centrefor Theoretical Physics in Trieste, Italy
%%}\\
%{\small(Preliminary version)}
}
\author{ 
%{\small (Preliminary version)}\\
%L. Chierchia \& 
C. E. Koudjinan \\
\vspace{-0.25truecm}
{\footnotesize Dipartimento di Matematica,  Universit\`a  ``Roma Tre"}
%\\ 
%\vspace{-0.25truecm}
%{\footnotesize Largo S. L. Murialdo 1, I-00146 Roma (Italy) }
\\
\vspace{-0.25truecm}
{\footnotesize ckoudjinan@mat.uniroma3.it}
\\ 
\vspace{-0.25truecm}
}
\maketitle
\begin{abstract}
%%We prove the persistence of a Cantor--family of KAM tori  of measure $O(\sqrt{\vae}^{1-2d/l}\qquad\vae^{1/2-d/l})$ for any Hamiltonian $H=K(y)+   P(y,x)$, with $K,P\in C^l(\mathscr D\times \tn)$ and $\mathscr D\subset \rn$ a bounded domain, provided $K$ is non--degenerate, $l>2d$, and  $\vae\coloneqq \|P\|_{C^l(\mathscr D\times \tn)}$ is sufficiently small \textcolor{red}{ or (   
%% provided $\vae^{?????}<c\; \a^2$ for some small constant $c=c(d,\t,K)>0$)}. %Moreover, the persistent KAM tori leave in a $\vae/\a$--neighborhood of their corresponding unperturbed tori. 
%%Furthermore, the regularity assumption can not be improved as shown in \cite{cheng2013destruction, wang2014total}.\\

%\frac{1}{2}-\frac{\n}{l}
\noi
{Given $l>2\n>2d\ge 4$, we prove the persistence of a Cantor--family of KAM tori  of measure $O(\vae^{1/2-\n/l})$ for any non--degenerate nearly  integrable Hamiltonian system of class $C^l(\mathscr D\times\tn)$, where $\n-1$ is the Diophantine power of the frequencies of the persitent KAM tori and $\mathscr D\subset \rn$ is a bounded domain, provided that the size $\vae$ of the perturbation is sufficiently small. This extends a result by D. Salamon in \cite{salamon2004kolmogorov} according to which we do have the persistence of a single KAM torus  %under the same regularity assumption;.
in the same framework. Moreover, it is well--known that, for the persistence of a single torus, the regularity assumption %is sharp as shown in \cite{}.
can not be improved.}
\end{abstract}
\tableofcontents
%\chapter{Quantitative KAM normal forms}
\section{Introduction}
KAM Theory asserts that, for sufficiently regular non--degenerate nearly integrable Hamiltonian systems, a Cantor--like family of KAM tori of the unperturbed part survive any perturbation, being only slightly deformed, provided the perturbation is small enough. Moreover, the family of KAM tori of the perturbed system is of positive Lebesgue measure and tends to fill up the phase space as the perturbation tends to zero. A natural question is:
\qst{qest1}
In a fixed degrees of freedom $d$, how regular has to be the integrable  Hamiltonian and the perturbation in order to get KAM tori?
\eqst
It was Arnold \cite{ARV63}, inspired by the breakthrough of Kolmogorov \cite{Ko1954},  who first proved the persistence of positive measure set of KAM tori of a real-analytic integrable Hamiltonian under a real-analytic perturbation, provided the latter is small enough.
In 1962, J. Moser\cite{moser1962new, moser1962invariant} proved in the framework of area--preserving twist mappings of an annulus, the persistence of invariant curves of integrable analytic systems under $C^k$ perturbation, but for $k$ very high ($k=333$); which, later on, was brought down by H. R\"ussmann \cite{russmann1970uber} to $5$ and finally to the optimal value 3 by M.R. Herman \cite{herman1983courbes}.%, which is  close to the optimal value 4. %\textcolor{red}{4} $3+\iota$, $\iota>0$. 

\noi
 Moser \cite{moser1969construction} proved the continuation of a single torus of an integrable real--analytic Hamiltonian under a perturbation of class $C^{l+2}$, with $l>2d$. Then P\"oschel \cite{poschel1980invariante,poschel1982integrability}, following an idea  due to Moser, showed that a Cantor--like family of KAM tori, of positive measure, of a non-degenerate integrable real--analytic Hamiltonian survive any sufficiently small perturbation of class $C^k$, provided $k>3d-1$, and also showed that, for the persistence of a single torus of an integrable real--analytic Hamiltonian, it is sufficient to require the perturbation to be of class $C^l$, provided $l>2d$. Later, refining this idea of Moser, D. Salamon \cite{salamon2004kolmogorov} showed that, for the persistence of single torus, it sufficient that both of the integrable and perturbed part are of class $C^l$, with $l>2d$. And, regarding the continuattion of a single torus of the integrable system, the regularity assumption $l>2d$ turns out to be also sharp (see \eg \cite{herman1983courbes,cheng2013destruction}). %indeed, Cheng and Wang \cite{cheng2013destruction} proved that any invariant Lagrangian torus with a given rotation vector of the integrable Hamiltonian $K_0=\sum_{j=1}^d y_j^2$ can be destroyed by arbitrarily $C^{2d-\d}$--small perturbations. 
 Indeed, M. Herman \cite{herman1983courbes} gave a counterexample of an exact area--preserving twist mappings of an annulus for which a given invariant curve  can be destroyed by arbirtrarly small perturbation of class $C^{3-\iota}$. It corresponds in the Hamiltonian context to $d=2$ and $l=4-\iota$. In \cite{herman1983courbes}, he also provided a counterexample of an exact area--preserving twist mappings of an annulus  of class $C^{2-\iota}$ without any invariant curve, corresponding in our context to $d=2$ and $l=3-\iota$.
Then, it has been widespread that
\conj{conjec1}
In $d$--degrees of freedom, a small perturbation of class $C^l$ of a non--degenerate integrable Hamiltonian which is also of class $C^l$, exhibits a positive measure set of KAM tori \ifof $l>2d$. 
\econj
Albrecht has proven in \cite{albrecht2007existence} the persistence of KAM tori of a non--degenerate real--analytic integrable system under small enough perturbations of class $C^{2d}$, provided that the moduli of continuity of the $2d$--th partial derivatives of the
perturbation satisfy some integral condition, weaker that the H\"older continuity condition. Yet, the KAM tori of the perturbed system form a zero measure set.

\noi
Given $\a>0$, $\t>0$, a vector $\o\in\rn$ is called $(\a,\t)$--Diophantine if\footnote{The positive constant $\t$ (resp. $\a$) is then called a Diophantine power (resp. Diophantine constant) of $\o$.}
$$
|\o\cdot k|\geq \frac{\a}{|k|_1^\t}, \ \ \forall\ k\in \zn\setminus\{0\}.
$$
In this paper, we prove the ``if'' part of the Conjecture~\ref{conjec1} \ie, roughly speaking:
\thm{teo00}
Consider a Hamiltonian of the form $H(y,x)=K(y)+P(y,x)$ where $K,P\in C^l(\mathscr D\times\tn)$ and $\mathscr D\subset \rn$ is a non--empty and bounded domain.\footnote{A domain is an open and connected set.}  If $K$ is non--degenerate and $l>2\n>2d$ then, %a positive measure set, say $\mathscr K$, formed by 
all the KAM tori of the integrable system $K$ whose frequency are $(\a,\t)$--Diophantine, %\footnote{See \equ{dio} hereafter.}
 with $\a\simeq \vae^{1/2-\n/l}$ and $\t\coloneqq \n-1$, do survive, being only slightly deformed, where $\vae$ is the $C^l$--norm of the perturbation $P$. % up to a small deformation. 
Moreover, letting $\mathscr K$ be the corresponding family of KAM tori of $H$, we have %for any bounded domain $\mathscr D\subset \rn$, 
$\meas(\mathscr D\times\tn\setminus \mathscr K)=O(\vae^{1/2-\n/l})$.
\ethm
To our best knowledge, the best result in this direction is due to A. Bounemoura and consigned in his nice paper \cite{bounemoura2018positive}, where he proved the persistence of positive measure set filled by the KAM tori of $H$ 
 %a family of KAM tori living out of a set of measure  $O(\sqrt{\vae})$  
 under the assumptions $K\in C^{l+2}$ and $P\in C^l$, with $l>2d$. Bounemoura also proved that the region free of KAM tori is of measure  $O(\sqrt{\vae})$ which turns out to be sharp (see \eg \cite{koudjinan2019quantitative,Chierrchia2019StructKolm}).
 %However this latter is not very clear to us. Indeed, in his proof, the Diophantine constant %\footnote{See \equ{dio} hereafter.} $\a$ (corresponding to $\g$ in \cite{bounemoura2018positive}) has been rescaled to one and this does not allow to keep track of the power of $\a$ relatively to $\vae$ which is crucial for the measure estimate. The point is that many other parameters of the KAM scheme, such as the analyticity domain of each of the real--analytic approximations of the perturbation, do depend upon $\vae$ and that need to be taken into account in the smallness condition and, in particular, in the measure estimate.  
 
 \noi
  In the present paper, under the sharper assumption $K,P\in C^l$, %and reasonable carefullness, 
we show that  the measure estimate of the region free of KAM tori %we are able to get 
is of $O(\vae^{1/2-\n/l})$, which, netherless, yields in the limit $l\rightarrow \infty$ the optimal bound in the real--analytic case \ie $O(\sqrt{\vae})$ ( see \eg \cite{koudjinan2019quantitative,Chierrchia2019StructKolm}).

\noi
The proof shares two main features with \cite{bounemoura2018positive}.  Firstly, our proof uses also a quantitative  approximation method of smooth functions by analytic functions introduced by Moser; however, here we have to approximate not only the pertubed part, but also the integrable part  at each step of the KAM scheme as, unlike \cite{bounemoura2018positive}, we do not linearize the integrable part. Secondly, we also use the refined approximation  given in \cite[Theorem~7.2, page~134]{russmann2001invariant} instead of truncating the Fourier expansion of the perturbation at each step of the KAM scheme. But, unlike \cite{bounemoura2018positive}, in this paper we use a KAM scheme \`a la Arnold.\footnote{Usually, in the literature, Moser's idea is combined with his own KAM scheme (like  \cite{poschel1980invariante,poschel1982integrability,bounemoura2018positive} ) or with Kolmogorov scheme (like \cite{salamon2004kolmogorov, CL12}).} The strategy is to prove a general quantitative  KAM Step for real--analytic perturbation of non--degenerate real--analytic integrable Hamiltonian systems (see Lemma~\ref{lem:1Extv5Simpl01}). Then, one approximates, in a quantitative manner, both the integrable and perturbed part by a sequence of real--analytic functions on complex strips of widths decreasing to zero (see Lemma~\ref{JacMosZen}), yielding a suitable real--analytic approximation of the perturbed Hamiltonian, to each of which we apply the KAM Step. Then, one proves that indeed the procedure converges.

\section{Notation\label{parassnot}}
%Fix $d\ge 2$. %Let $\O\subset \rn$ be non-empty and bounded domain with piecewise smooth boundary and 
%We shall use the following standard notation.
%
%In the following, we let
%$r,\,s,\,\a,\,\t,\,N>0;\, n,p\in \natural\coloneqq \{1,2,3\cdots\}$; $ y_0\in \complex^n$. 

\begin{itemize}

\item[\tiny $\bullet$] For $d\in\natural \coloneqq \{1,2,3,...\}$ and  $x,y\in \complex^d$, we let
$x\cdot y\coloneqq x_1 \bar y_1 +\cdots+x_d \bar y_d$ be the standard inner product;
$|x|_1\coloneqq\dst\sum_{j=1}^d |x_j| $ be the $1$--norm,
and  $|x|\coloneqq \dst\max_{1\le j\le n}|x_j|$ be the sup--norm.

\item[\tiny $\bullet$] $\tn\coloneqq \rn/2\p\zn$ is the  $d$--dimensional (flat) torus.

\item[\tiny $\bullet$] For $\a>0$, $\t\ge d-1\ge 1$, 
\beq{dio}\D_\a^\t\coloneqq \left\{\o\in \rn:|\o\cdot k|\geq \frac{\a}{|k|_1^\t}, \ \ \forall\ 0\not=k\in \zn\right\},
\eeq
is the set of $(\a,\t)$--Diophantine numbers in $\real^d$.

\item[\tiny $\bullet$] We denote by $\meas$, the Lebesgue (outer) measure on $\rn$;

\item[\tiny $\bullet$] Given $l\in \real$, we shall denote its integer part by $[l]$ and its fractional part by $\{l\}$;

\item[\tiny $\bullet$] For %$\mathscr D\subseteq \rn$ and 
$l>0$, $A$ an open subset of $\rn$ or of $\rn\times \tn$, we denote by $C^l(A)$ the set of continuously differentiable functions $f$ on $A$ up to the order $[l]$ such that $f^{[l]}$ is H\"older--continuous with exponent $\{l\}$ and with finite $C^l$--norm define by:
\begin{align*}
&\|f\|_{C^l(A)}\coloneqq \max\{\|f\|_{C^{[l]}(A)}\;,\,\|f^l\|_{C^{\{l\}}(A)}\}\;, \quad\|f\|_{C^{[l]}(A)}\coloneqq\max_{\substack{k\in\nn\\0\le |k|_1\le [l]}}\sup_{A}|\dpr^{k}_y f| \;,\\
& \|f^l\|_{C^{\{l\}}(A)}\coloneqq \max_{\substack{k\in\nn\\ |k|_1= [l]}}\;\sup_{\substack{y_1,y_2\in A\\ 0<|y_1-y_2|<1}}|\dpr^{k}_yf(y_1)-\dpr^{k}_yf(y_2)|/|y_1-y_2|^{\{l\}}\;.
  %&\hspace{7cm} z_1,z_2\in\rn\times\tn\;,\, 0<|y_1-y_2|<1\bigg\}
\end{align*}
When $A=\rn$ or $A=\rn\times \tn$, we will simply write $\|f\|_{C^l}$ for $\|f\|_{C^l(A)}$.

\item[\tiny $\bullet$] For $l>0$, $A$ any subset of $\rn$, we denote by $C_W^l(A)$, the set of functions of class $C^l$ on $A$ in the sense of Whitney.\footnote{We refer the reader for instance to \cite[Appendix~E, page~207]{koudjinan2019quantitative} for details.}

\item[\tiny $\bullet$] For $r,s>0$, $y_0\in\complex^d$, $\emptyset\neq\mathscr D\subseteq \cn$, we denote:
\beqano
\dst\torus^d_s &\coloneqq  &\left\{x\in \cn: |\Im x|<s\right\}/2\p\zn\,,\\
B_r(y_0)&\coloneqq  &\left\{y\in \real^d: |y-y_0|<r\right\}\,,\qquad (y_0\in\real^d)\,,\\
D_r(y_0)&\coloneqq  &\left\{y\in \cn: |y-y_0|<r\right\}\,,\\
D_{r,s}(y_0)&\coloneqq  & D_r(y_0)\times \torus^d_s\,,\\
D_{r,s}(\mathscr D)&\coloneqq & \dst\bigcup_{y_0\in\mathscr D}D_{\mathsf{r},s}(y_0)\,.
\eeqano 

\item[\tiny $\bullet$] If  $\uno_d\coloneqq \diag(1)$ is the unit $(d\times d)$ matrix, we denote the standard symplectic matrix by 
$$\mathbb{J}\coloneqq  \begin{pmatrix}0 & -\uno_d\\
\uno_d & 0\end{pmatrix}\,.
$$ 

\item[\tiny $\bullet$]
For $\mathscr D\subset \cn$, $\mathcal{A}_{r,s}(\mathscr D)$ denotes the Banach space of real--analytic functions with bounded holomorphic extensions to $ D_{r,s}(\mathscr D)$, with norm 
$$
\|\cdot\|_{r,s,\mathscr D}\coloneqq \dst\sup_{ D_{r,s}(\mathscr D)}|\cdot|\,.
$$
\item[\tiny $\bullet$] We equip 
$\cn\times\cn$   with the canonical symplectic form 
$$\varpi\coloneqq dy\wedge dx=dy_1\wedge dx_1+\cdots+dy_d\wedge dx_d\ ,
$$
and denote by $\phi_H^t$ the associated Hamiltonian flow governed by the Hamiltonian $H(y,x)$, $y,x\in\complex^d$.
\item[\tiny $\bullet$]
 $\pi_1\colon \cn\times\cn\ni(y,x)\longmapsto y$ is the projection on the first $d$--components and, $\pi_2\colon \cn\times\cn\ni(y,x)\longmapsto x$ is the projection on the last $d$--components.
\item[\tiny $\bullet$]
Given a linear operator $\mathcal{L}$ from the normed space $(V_1,\|\cdot\|_1)$ into the normed space  $(V_2,\|\cdot\|_2)$, its ``operator--norm'' is given by
\[\|\mathcal{L}\|\coloneqq \sup_{x\in V_1\setminus\{0\}}\frac{\|\mathcal{L}x\|_2}{\|x\|_1},\quad \mbox{so that}\quad \|\mathcal{L}x\|_2\le \|\mathcal{L}\|\, \|x\|_1\quad \mbox{for any}\quad x\in V_1.\]
\item[\tiny $\bullet$]
Given $\o\in \rn$, the directional derivative of a $C^1$ function $f$ with respect to $\o$ is given by
 
\[D_\o f\coloneqq \o\cdot  f_x=\dst\sum_{j=1}^d \o_j \dst f_{{x}_j}\,.\]
\item[\tiny $\bullet$]
If  $f$ is a (smooth or analytic)  function on $\torus^d$, its Fourier expansion is given by    \[f=\dst\sum_{k\in \zn}f_k \ex^{\mathbf{i} k\cdot x}\,,
\qquad f_k:=\dst\frac{1}{(2\pi)^d}\dst\int_{\tn}f(x) \ex^{-\mathbf{i} k\cdot x}\, d x\,,\]
(where, as usual, $\ex\coloneqq \exp(1)$ denotes the Neper number and $\mathbf{i}$ the 
imaginary unit). We also set:
\[\average{f}\coloneqq f_0=\dst\frac{1}{(2\pi)^d}\dst\int_{\tn}f(x)\, d x\,.%\qquad  T_N f:=\dst\sum_{|k|_1\leq N}f_k \ex^{ik\cdot x},\, N>0. 
\]

\end{itemize}

\section{Assumptions\label{settings}}

\begin{itemize}
\item[\tiny $\bigstar$] Let $l>2\n\coloneqq 2(\t+1)>2d \ge 4$. and $\mathscr D\subset\rn$ be a non--empty, bounded domain. 
%As we are interested in a weak regularity assumption, we may consider $l$ very close to $2\n$ and, therefore,  assume\footnote{Notice that, the assumption $l<2\n+1$ is made just for simplicity and may be removed  by changing suitably the definition of $\s$ in \equ{SigmRoBet0} (see {\bf (ii)} of Remark~\ref{DstD0Fi0}).}
%\beq{lNu}
%2\n<l<2\n+1.
%\eeq
\item[\tiny $\bigstar$] On the phase space $\mathscr D\times\torus^d$, consider the Hamiltonian
\[{\mathrm{H}}(y,x)\coloneqq \mathrm{K}(y)+  {\mathrm{P}}(y,x),\]
where  $\mathrm{K}, {\mathrm{P}}\in C^l(\mathscr D\times\tn)$ are  given functions with finite $l$--norms $\|\mathrm{K}\|_{C^l(\mathscr D)}$ and $\vae\coloneqq \|{\mathrm{P}}\|_{C^l(\mathscr D\times \tn)}$.% and $\vae$ a positive real parameter.
\item[\tiny $\bigstar$] Assume that $\mathrm{K}_y$ is locally--uniformly invertible; namely that $\det \mathrm{K}_{yy}(y)\neq 0$ for all $y\in \mathscr D$ and 
$$
%\mathsf{T}\coloneqq
%\|T\|_{\mathscr D}
\mathsf{T}\coloneqq \|T\|_{C^0(\mathscr D)}<\infty,\qquad T(y)\coloneqq \mathrm{K}_{yy}(y)^{-1}.
$$
Set\footnote{Indeed, $\eta\ge \|T(y_0)\|\|\mathrm{K}_{yy}(y_0)\|=\|T(y_0)\|\|T(y_0)^{-1}\|\ge 1$, for any $y_0\in \mathscr D$.} 
$$\mathsf{K}\coloneqq\max\left\{1,\|\mathrm{K}\|_{C^l(\mathscr D)}\right\},\quad \eta\coloneqq \mathsf{T}\mathsf{K}\ge 1.
%\quad \varrho\coloneqq \dst\inf_{y\in \mathscr{D}}|\det K_{yy}(y)|> 0,\quad \vth\coloneqq \frac{\mathsf{K}^d}{\varrho}\ge 1.
$$
\item[\tiny $\bigstar$] Let $\a\in (0,1)$ and set
$$
\a_*\coloneqq \a^{\frac{1}{l-2\n}}
%\left(\frac{\vae}{\mathsf{K}}\right)^{\frac{1}{2}}\exp\left(\n C_5\eta^{\su\n}\right)
\;,\qquad\mathscr D' \coloneqq \left\{y\in\mathscr D:  B_{\a_*}(y)\subseteq\mathscr D\right\}%\quad\mbox{and}\quad
$$
and
$$\mathscr D_{\a} \coloneqq \left\{y\in\mathscr D': \ \mathrm{K}_y(y)\in \D_\a^\t\right\}\,.
$$
%%small enough so that the set
%%$
%%\mathscr D'\coloneqq \{y\in \mathscr D\ :\ {B_{2\a}(y)}\subseteq \mathscr D\}
%%$
%%is of non--empty interior. Set
%%$$
%%\D\coloneqq \left\{y\in \mathscr D'\ :\ K_y(y)\in \D_\a^\t \right\}.
%%$$
\item[\tiny $\bigstar$]
%%Let introduce the rescaled smallness parameter $\r$:
%%$$
%%\r\coloneqq \frac{2\|\mathrm{K}\|_{C^l}\;\|{\mathrm{P}}\|_{C^l}}{\a^2\s^{3\n}}=\frac{2\mathsf{K}\;\vae}{\a^2\s^{2\n}}\;,
%%$$
%%where
%%$$
%%\textcolor{red}{\s\coloneqq \left(\frac{1}{\eta^{l\n+6}}\frac{\vae^2}{\a^2}\right)^{1/(l+4\n)}}.
%%$$
Finally, set 
\beq{SigmRoBet0}
\s\coloneqq \left(\frac{\vae^{3/2}}{\eta^{2l/\n}\a\sqrt{\mathsf{K}}}\right)^{1/(l+\n)}\;,\qquad \r\coloneqq \frac{2C_1\mathsf{K}\;\vae}{\a^{2}\s^{2\n}},\qquad
\b_0\coloneqq \min\left\{\frac{l}{2\n}-1+\su\n\;,\;2\right\}\;,
\eeq
for some suitable constant $C_1=C_1(d,l)>1$.\footnote{$C_!$ is actually the constant apperaring in Lemma~\ref{JacMosZen}.}
%and introduce the smallness parameter
\end{itemize}
\section{Theorem\label{SecStatem}}
Under the notations and assumptions of \S~\ref{parassnot} and \ref{settings}, the following Theorem holds.
\thm{Extteo4v2}\ \\ %{Arnold \cite{ARV63}}
{\bf Part I: }
There exists a positive constant $\mathsf{c}=\mathsf{c}(d,\t,l)<1$  such that, if%\footnote{Observe that $\equ{smcEAr0v2} \Longleftrightarrow \r\le c$.}
\beq{smcEAr0v2}
% \vae\coloneqq \|{\mathrm{P}}\|_{C^l(\mathscr D\times \tn)}\le \mathsf{c}\;\frac{\ex^{-\eta}\a^2\s_0^{3\n}}{\|\mathrm{K}\|_{C^l(\mathscr D)}}
\a\le \mathsf{c}\;\mathsf{K}\;,\qquad \mbox{and}\qquad\vae\le \mathsf{c}\;\mathsf{K}^{-\frac{l+2\n}{l-2\n}}\eta^{-a}\a^{\frac{2l}{l-2\n}}\;,
\eeq
where 
$a\coloneqq (l-2\n)^{-1}\max\{(6+2l\n^{-1})(l+\n)-2l(l-\n),\; 2l(l+3\n)\n^{-1}\},
$ 
then, the following holds.
%% with
%%$$
%%\textcolor{red}{r_*		   \coloneqq \frac{\s^\n}{\mathsf{C}_{9}}\left(\frac{\s}{\eta}\right)^{\frac{5}{4}}\frac{\a}{\mathsf{K}}\;.}
%%$$
 There exist a Cantor--like set $\mathscr D_*\subset {\mathscr D}$, an embedding $\phi_*=(v_*,u_*)\colon \mathscr D_*\times \tn\to \mathscr K\coloneqq \phi_*(\mathscr D_*\times \tn)\subset \mathscr D\times \tn$ of class $C_W^{\b}(\mathscr D_*\times \tn)$ such that the map $\x\longmapsto \phi_*(y_*,x)$ is of class $C^{\b\n}(\tn)$, for any given $y_*\in\mathscr D_*$ (for any $\n^{-1}<\b<\b_0$), a function $\mathrm{K}_*\in C_W^2(\mathscr D_*,\real)$, satisfying 
\begin{align}
{\mathrm{H}}\circ \phi_*(y_*,x)&=\mathrm{K}_*(y_*),\qquad \forall\;(y_*,x)\in \mathscr D_*\times\tn. \label{conjCaneq0v2}
\end{align}

\noi
Moreover, the map $G^*\coloneqq (\dpr_{y_*}\mathrm{K}_*)^{-1}\circ\dpr_{y}\mathrm{K} \colon \mathscr D_{\a}\longrightarrow\mathscr D_*$ is well-defined and is a lipeomorphism onto $\mathscr D_*$, $B_{\a_*/2}(\mathscr D_*)\subseteq \mathscr D$, and $\mathscr K$ is foliated by KAM tori of $H$, each of which is a graph of a map of class $C^\n(\tn)$.\footnote{See {\bf (i)} in Remark~\ref{DstD0Fi0} below.} Furthermore,
\begin{align}
&\|G^*-\id\|_{\mathscr D_{\a}}
\le %\textcolor{red}
{\vae^{\frac{3\t}{2(l+\n)}}\a^{\frac{l+1}{l+\n}}\mathsf{K}^{\frac{-\t}{2(l+\n)}}\eta^{-1-\frac{2l\t}{\n(l+\n)}}}\,,\qquad \|G^*-\id\|_{L,\mathscr D_{\a}}<1/2\,,\label{NormGstrThtv2}\\
\ \nonumber\\
&\sup_{\mathscr D_*\times\tn}\max\{|\mathsf{W}(\phi_*-\id)|\;,\|\pi_2(\dpr_x \phi_*-\uno_d)\|\}\le 8\eta^{-2}(\log\r^{-1})^{-2\n}<1
\;, \label{estArnTrExtv2}
\end{align}
where %$\r\coloneqq (\vae^{l-2\n}\a^{-2\n}\mathsf{K}^{l+2\n}\eta^{4l})^{1/(l+\n)}$ and
$
\mathsf{W}   \coloneqq \diag({\mathsf{K}}({{\a\s^\n}})^{-1}\uno_d,\s^{-1}\uno_d).
$
\ \\

\noi
{\bf Part II: } 
Assume furthermore that the boundary $\dpr \mathscr D$ of $\mathscr D$ is a smooth hypersurface of $\rn$ and
\beq{EqStnAlf}
0< \a\le \min\left\{\frac{R(\mathscr D)}{6}\,,\,\su 2\minfoc(\dpr \mathscr D)\right\}\;,
\eeq
where $\minfoc(\dpr \mathscr D)$ denotes the minimal focal distance of $\dpr \mathscr D$ and\footnote{Observe that the condition $\a\le R(\mathscr D)/6$ ensures that the interior of $\mathscr D'$ is non--empty.}
$$
R(\mathscr D)\coloneqq \sup\{R>0: B_R(y)\subseteq \mathscr D\;, \mbox{ for some } y\in\mathscr D\}\;.
$$
 Then, the following measure estimate holds:
\beq{MesChPin}
\meas(\mathscr D\times\tn\setminus \mathscr K)\le (3\pi)^d \bigg(2\mathcal{H}^{d-1}(\dpr \mathscr D)\;\wh \vae+ C\;{\wh \vae}^2+\meas(\mathscr D'\setminus \mathscr D_{\a})\bigg)\,,%\marginnote{to be rechecked!!! $\a^?$}
\eeq
where\footnote{We refer the reader to \cite{Chierrchia2019StructKolm, koudjinan2019quantitative} for more details.
} $\mathbf{R}^{\dpr\mathscr D}$ denotes the curvature tensor of $\dpr\mathscr D$, $\mathbf{k}_{2j}(\mathbf{R}^{\dpr\mathscr D})$, the $(2j)$--th integrated mean curvature of $\dpr\mathscr D$ in $\rn$,
$$
\wh \vae\coloneqq \max\left\{{\vae^{\frac{3\t}{2(l+\n)}}\a^{\frac{l+1}{l+\n}}\mathsf{K}^{\frac{-\t}{2(l+\n)}}\eta^{-1-\frac{2l\t}{\n(l+\n)}}}\,,\ \a_* \right\}\,,
$$
 and 
$$
C=C(d,\t, l, \vae, \a,\mathsf{T}, \mathsf{K}, \mathbf{R}^{\dpr\mathscr D})\coloneqq 2\dst\sum_{j=1}^{\left[\frac{d-1}{2}\right]}\frac{{\wh \vae}^{2j-1}\mathbf{k}_{2j}(\mathbf{R}^{\dpr\mathscr D})\;}{1\cdot3\cdots (2j+1)} \;.
$$
\ethm
\rem{DstD0Fi0}
{\bf (i)} 
By definition,
\beq{conjCaneq00v2}
\dpr_{y_*}\mathrm{K}_*\circ G^*=\dpr_{y}\mathrm{K}  \qquad \mbox{on} \quad \mathscr D_{\a}\;.
\eeq
Now, from \equ{conjCaneq00v2} and \equ{conjCaneq0v2}, one deduces that the embedded $d$--tori
\beq{KronTorArnExt}
\mathcal{T}_{\o_*,\vae}\coloneqq \phi_*\left(y_*,\tn\right),\qquad y_*\in \mathscr D_*\,,\quad \o_*\coloneqq \dpr_{y_*}K_*(y_*)\in \D_\a^\t\,,
\eeq
are non--degenrate, invariant, Lagrangian Kronecker tori of class $C_W^{{\b\n}}$ (for any $\n^{-1}<\b<\b_0$) for $H$, \ie KAM tori, with Diophantine frequency $\o_*$ \ie
\beq{KronTorArnIEExt}
\phi^t_H\circ \phi_*(y_*,x)=\phi_*(y_*,x+\o_* t)\,, \qquad \forall\; x\in\tn.
\eeq
%and, by \equ{estArnTrExtv2}, they are all graphs over the ``angle'' variables of a function of class $C_W^{??}$. Thus, each $\mathcal{T}_{\o_*,\vae}$ is a KAM torus of class $C_W^{???}$ of $H$.
Indeed, as each $\phi^j$ is symplectic, we have
\beq{ConjHamSym}
\phi^t_{\mathcal{H}_{j-1}}\circ\phi^{j}=\phi^{j}\circ \phi^t_{\mathcal{H}_{j-1}\circ\phi^{j}}.
\eeq
Now, pick $y_*\in\mathscr D_*$ and $y_j\in\mathscr D_j$ converging to $y_*$. Letting $\o_*\coloneqq \dpr_{y_*}K_*(y_*)$, we have
\beq{ApproTori}
\phi^t_{\mathcal{H}_{j-1}\circ\phi^{j}}(y_j,x)\eqby{KAMToriH} (y_j,x+t\o_*)+O(r_j^{-1}\|P_j\|_{r_j,s_j,\mathscr D_j}+|y_j-y_*|)\;,\qquad \lim_{j\rightarrow \infty}r_j^{-1}\|P_j\|_{r_j,s_j,\mathscr D_j}\eqby{estfin2Ext01v501}0.
\eeq
 Then, recalling that $\mathcal{H}_j$ converges uniformly to $H$ on $\rn\times\tn$, we have, for any $x\in\tn$,
\begin{align*}
\phi_H^t\circ \phi_*(y_*,x)&=\lim_{j\rightarrow \infty}\phi_{\mathcal{H}_{j-1}}^t\circ \phi_j(y_j,x)\\
		&\eqby{ConjHamSym} \lim_{j\rightarrow \infty}\phi^{j}\circ \phi^t_{\mathcal{H}_{j-1}\circ\phi^{j}}(y_j,x)\\
		&\eqby{ApproTori} \phi_*(y_*,x+t\o_*),
\end{align*}
and \equ{KronTorArnIEExt} is proven. In particular, each torus $\mathcal{T}_{\o_*,\vae}$ is Lagrangian by Lemma~\ref{LagTor}.\footnote{An alternative proof is as follows. Let $v_j\coloneqq \pi_1\phi_j$ and $u_j\coloneqq \pi_2\phi_j$. Then, as $\phi_j$ is symplectic, we have $(\dpr_x u_j)^T \dpr_x v_j=(\dpr_x v_j)^T \dpr_x u_j$, and letting $j\to \infty$, we get $(\dpr_x u_*)^T \dpr_x v_*=(\dpr_x v_*)^T \dpr_x u_*$ \ie the tori $\mathcal{T}_{\o_*,\vae}$ are Lagrangian, where $A^T$ denotes the transpose of the matrix $A$.} \qed

%%\noi
%%{\bf (ii)} Unlikely, the smallness condition \equ{smcEAr0v2} we obtain is not sharp. Indeed,
%%$
%% a= {(l-2\n)}/{(l+\n)}<1,
%% $
%% which is not what one expects; in fact, one expects the power of $\vae$ to be 1.

\noi
{\bf (ii)} Choosing %$\vae\simeq \a^{{2l}/{(l-2\n)}}$ 
$\a\simeq \vae^{1/2-\n/l}$ in \equ{smcEAr0v2}, we get $\wh \vae=O(\vae^{1/2-1/l})$ and therefore, plugging them into  \equ{MesChPin}, we obtain 
\beq{SimplMeas}
\meas(\mathscr D\setminus\mathscr K)=O(\vae^{\su2-\frac{\n}{l}}),
\eeq
which agrees for $l\rightarrow\infty$ with the sharp measure of KAM tori for smooth Hamiltonian systems \ie $O(\sqrt{\vae})$. Moreover, \equ{NormGstrThtv2} yields the sharp bound $O(\vae^{\su2-\frac{\n}{l}})$ on the displacement of each persistent invariant torus from the corresponding unperturbed one \ie 
$$
%\sup_{\mathscr D_\a\times \tn}|v_*-\id|
\|G^*-\id\|_{\mathscr D_{\a}}=O(\vae^{\su2-\frac{\n}{l}}).
$$
The relation in  \equ{estArnTrExtv2} yields oscillations of $O(\sqrt{\vae})$ for each perturbed torus:
$$
\sup_{y_*\in \mathscr D_*}\sup_{x,x'\in \tn}|v_*(y_*,x)-v_*(y_*,x')|=O(\sqrt{\vae}),
$$
which is sharp (see \eg \cite{chierchia2019ArnoldKAM}). It is worth mentioning that, in order to get \equ{SimplMeas}, the smoothness assumption on the boundary of the domain can be removed using a different argument. The argument consists in slicing the domain into small pieces, then construct in each of those pieces a family of KAM tori and estimate their respective relative measures, and finally some them all up (see \cite{Chierrchia2019StructKolm, koudjinan2019quantitative} for more details).
%% \noi
%% {\bf (iii)} Actually, $G^*\in C^{1+\m}_W(\mathscr D_*)$, for any $0<\m<l/\n-2$.
\erem

%\section{Proof}
\section{Proof of Theorem~\ref{Extteo4v2}}

\subsection{General step of the KAM scheme}
\lem{lem:1Extv5Simpl01}
Let $r>0,\; 0<\bar{\s}\le 1,\; 0< 2\s< s\le 1$, $\mathscr D_\sharp\subset\rn$ be a non--empty, bounded domain. Consider the Hamiltonian
$$
H(y,x)\coloneqq K(y)+ P(y,x)\;,
$$
where $K,P\in \mathcal{A}_{r,s}(\mathscr D_\sharp)$. % $K,P$ are real--analytic functions with bounded holomorphic extensions to $D_{r,s}(\mathscr D_\sharp)$.\\
Assume that%\footnote{In the sequel, $K$ and $P$  stand for  generic real analytic Hamiltonians which, later on, will respectively play the roles of $K_j$ and $P_j$,  and $y_0,\,r$, the roles of $y_j,\,r_j$ in the iterative step.} 
\beq{RecHypArnExtv501}
\begin{aligned}
%&r\le r_0\;, \qquad\qquad\qquad\qquad\qquad \mathscr D_\sharp\subset \mathscr D_r 
	%\left\{y_0\in\mathscr D:\ \dist(y_0,\dpr\mathscr D)\ge \frac{r}{32d}\ \mbox{ and }\ K_y(y_0)\in \D_\a^\t 
	%\right\}
%\;,\\ 
&%\dst\inf_{y\in\mathscr D}|\det K_{yy}(y)|\ge \d>\d_\infty%\left(1-\frac{\s}{6}\right)^{-d}
\det K_{yy}(y)\not=0
\;,\qquad\qquad\qquad T(y)\coloneqq K_{yy}(y)^{-1}\;,\quad \forall\;y\in\mathscr D_\sharp\;,\\
%&r_0|\o|,\, r_0\|K_y\|_{r,y_0},\, r_0^2\|K_{yy}\|_{r,y_0},\, |\o|^2\|T\|\le \su2\mathsf{E}\le \su2\mathsf{E}_\infty,\\ %\quad \a<1,\quad r<1,
&\|K_{yy}\|_{r,\mathscr D_\sharp}\le \mathsf{K}\;,\qquad\qquad\qquad\ \, \|T\|_{\mathscr D_\sharp}\le \mathsf{T}\;,\\
& \|P\|_{r,s,\mathscr D_\sharp}\le \vae \;,\qquad\qquad\qquad\,\,\, K_y(\mathscr D_\sharp)\subset \D^\t_\a\;. %\qquad\qquad r\le r_0 \;,
\end{aligned}
\eeq
\noi
%with $\mathscr K\subset \mathscr D$ and $T=K_{yy}^{-1}$.
Assume that
\beq{DefNArnExt1v501}
\s^{-\n}\frac{\vae}{{\a}r}\le \r\le \su4
%\textcolor{red}{\leby{cond1ExtExtv5} \su4 \qquad To\ BE\ ERASED} 
\qquad\mbox{and}\qquad r\le \frac{\a}{\mathsf{K}}\s^\n
\;.
\eeq
Let
\beq{DefNArn2v501}
\begin{aligned}
&\eta\coloneqq\mathsf{T}\mathsf{K}\;,\quad\l\coloneqq \log\r^{-1}\;,\quad
\k\coloneqq 6\s^{-1}\l\;, \quad \check{r}\le \frac{r}{32d\mathsf{T} \mathsf{K}}\;,\quad 
\bar{r}%\le \frac{r\s^\n}{2d\cdot 6^\n\l^\n\eta}\textcolor{red}
{\le \dst\min\left\{\frac{\a}{2d\mathsf{K}\k^{\n}}\,,\, \check{r} \right\}},\\
&%\bar{\s}\le \frac{\s}{\s_0},\quad
\tilde r\coloneqq \frac{\check{r}\bar{\s}}{16d\mathsf{T} \mathsf{K}}
\;,\quad \bar{s}\coloneqq s-\frac{2}{3}\s\;,\quad s'\coloneqq s-\s \;, \qquad \mathsf{L}\coloneqq C_0\;\frac{\eta\mathsf{T}\vae}{r\tilde{r}}
\;.
\end{aligned}
\eeq
%%
%%Then, there exists a generating function $(y',x)\mapsto y'\cdot x+g(y',x)$, with $g\in \mathcal{B}_{\bar r,\bar s}(\mathscr D_\sharp)$
%%%and a function $\wt K\in \mathcal{B}_{\bar r,\bar s}(\mathsf{y})$,
%% and satisfying the following inequalities:
%%\beq{Est1Lem1bExtv5}
%%\left\{
%%\begin{aligned}
%%&\|g_x\|_{\bar{r},\bar{s},\mathscr D_\sharp}\le  \mathsf{C}_1 \frac{M}{\a} \s^{-(\n+d)}\,,\\
%%& \|g_{y'}\|_{\bar{r},\bar{s},\mathscr D_\sharp},\, \|\dpr_{y'x}^2 g\|_{\bar{r},\bar{s},\mathscr D_\sharp}\le \ovl{\mathsf{L}}\,,\\%\mathsf{C}_0 \max(\a,r\mathsf{K})\frac{M}{\a^2 r}\s^{-(2\t+d+2)}\,,\\
%%%,\, |\o| |\wt y|, |\o|^2\|\wt T\|
%%%&\|\dpr_{y'x}^2 g\|_{\bar{r},\bar{s},y_0}\le\mathsf{C}_0 \frac{r_0\mathsf{K}M}{\a^2 r}\s^{-(2\t+d+2)}\,, \\
%%&\|\dpr_{y'}^2\wt K\|_{\bar{r},\mathscr D_\sharp}\le %\frac{4M}{r^2}\,. 
%%\mathsf{K}\mathsf{L}\,,
%%\end{aligned}
%%\right.
%%\eeq
%%where
%%$$\wt K(y')\coloneqq \average{P(y',\cdot)}\;.$$
%%
Assume:   
\beq{cond1ExtExtv501}
%\vae_*\le 2\min\left(\vae_\sharp\,,\, \frac{r\bar{r}\s^{2(\t+d+1)}}{16\mathsf{T}_\infty  M}\,,\,\frac{\a r_0}{\mathsf{C}_2 M}\s^{4d}\,,\,\frac{\mathsf{K} r\s}{12M}\right)
{\mathsf{L}}\le \frac{\bar{\s}}{3}
\ .
\eeq
Then, there exists %$\mathscr K_1\subset\mathbb B^d_{\bar{r}}(\mathscr K)$ and 
a  diffeomorphism $G\colon  D_{\tilde{r}}(\mathscr D_\sharp){\to}G( D_{\tilde{r}}(\mathscr D_\sharp))$, a symplectic change of coordinates
\beq{phiokExt0Ext01}
\phi'=\id+  \tilde{\phi}: D_{\bar{r}/2,s'}(\mathscr D_\sharp')\to D_{\bar{r}+r\s/3, \bar{s}}(\mathscr D_\sharp),%\quad{\rm with}
%\quad 
\eeq
%generated by a function $g\in \mathcal{A}_{\bar{r},\bar{s}}(\mathscr D_\sharp)$ 
such that 
\beq{HPhiH'01}
\left\{
\begin{aligned}
& H\circ \phi'\eqqcolon H'\eqqcolon K'+ P'\ ,\\%\quad K'= K'(y'),\\
& \dpr_{y'}K'\circ G=\dpr_{y}K,\quad \det \dpr_{y'}^2K'\circ G\neq 0\quad\mbox{ on } \mathscr D_\sharp\,,
\end{aligned}
\right.
\eeq
with $K'(y')\coloneqq K(y')+ \wt K(y')\coloneqq K(y')+  \average{P(y',\cdot)}$. Indeed, $G=(\dpr_{y'}K')^{-1}\circ K_y$.  Moreover, letting $\left(\dpr^2_{y'} K'(\mathsf{y}')\right)^{-1}\eqqcolon T(\mathsf{y}')+ \;\wt T(\mathsf{y}')$, $\mathsf{y}'\in G(\mathscr D_\sharp)$, the following estimates hold.

\beq{convEstExt01}
%\left\{
\begin{aligned}
&\|\dpr_{y'}^2\wt K\|_{r/2,\mathscr D_\sharp}\le %\frac{4M}{r^2}\,. 
\mathsf{K}\mathsf{L}\,,\qquad \|G-\id\|_{\tilde{r},\mathscr D_\sharp}\le \tilde{r}\mathsf{L}\;,\qquad \|\wt T\|_{\mathscr D_\sharp'}\le \mathsf{T}\mathsf{L}\,, \\
& \max\{\|\mathsf{W}\,\tilde \phi\|_{\bar{r}/2,s',\mathscr D_\sharp'}\;,\,\|\pi_2\dpr_{x'}\tilde \phi\|_{\bar{r}/2,s',\mathscr D_\sharp'}\}\le C_1\frac{\vae}{\a r\s^\n}\,,\qquad \|P'\|_{\bar{r}/2, s',\mathscr D_\sharp'} \le  C_1\;\r \;\vae\,, %\qquad\qquad\quad\qquad\quad\ \,\,\, B_{\bar{r}/2}(\mathscr D_\sharp')\subset B_{\bar{r}}(\mathscr D_\sharp)\subset\mathscr D\,,
%\dst\inf_{\mathscr D}|\det\dpr_{y'}^2K_1|\ge \d\left(1-\frac{\s}{6}\right)^d\,,%>\d_\infty\,,
\end{aligned}
%\right. 
\eeq
where
\begin{align*}
&\mathscr D_\sharp'\coloneqq G(\mathscr D_\sharp)\;, \quad \left(\dpr^2_{y'} K'(\mathsf{y}')\right)^{-1}\eqqcolon T\circ G^{-1}(\mathsf{y}')+ \;\wt T(\mathsf{y}')\,, \ \forall\; \mathsf{y}'\in \mathscr D_\sharp'\,, \\
& \mathsf{W}\coloneqq \diag(r^{-1}\uno_d,\s^{-1}\uno_d)\,.
\end{align*}
%$$
%\mathsf{W}\coloneqq \begin{pmatrix}
%\max\{\frac{\mathsf{K}}{{\a}}\;,\frac{1}r\}\;\uno_d & 0\\ \ \\
%0			& \uno_d 
%\end{pmatrix}\;.
%$$
\elem
\proof
%\textcolor{red}
{The proof follows essentially the same lines as the one of the KAM Step in \cite{koudjinan2019quantitative} (see also \cite{Chierrchia2019StructKolm}) modulo two changes:\\%, see also \cite{koudjinan2019quantitative}.
{\bf (i)} To construct the generating function, as in \cite{bounemoura2018positive}, we use the approximation  given in \cite[Theorem~7.2, page~134]{russmann2001invariant} instead of truncating the Fourier expansion of $P$.\\
{\bf (ii)} We use systematically the estimate in 2. of Lemma~\ref{IFTLem} to estimate the generating function as well as its derivatives.\\
Those two modifications improve a lot the KAM Step; in particular it yields the optimal power of the lost of regularity $\s$, which is crucial in the KAM Theory for finitely differentiable Hamiltonian systems, at least from the Moser's ``analyticing'' idea point of view.
} We refer the reader to Appendix~\ref{appC} for an outline of the proof.
\qed
\subsection{Characterization of smooth functions by mean of real--analytic functions}
The following two Lemmata, which will be needed from Lemma~\ref{IteKAM} on and  may be found in \cite{chierchia2003kam, salamon2004kolmogorov}.
 
\lemtwo{JacMosZen}{Jackson, Moser, Zehnder}
Given $l>0$, there exists $C_1=C_1(d,l)>0$ such that for any $f\in C^l(\rn)$ and for any $\mathsf{s}>0$, there exists a real--analytic function $f_\mathsf{s}\colon \mathcal{O}_s\coloneqq \{(y,x)\in \cn\times\cn\,:\,|\Im(y,x)|<\mathsf{s}\}\to \complex$ satisfying the following:
\beq{AnalAppro}
\sup_{\mathcal{O}_\mathsf{s}}|f_\mathsf{s}|\le C_1\|f\|_{C^0},\quad
% \sup_{\mathcal{O}_{\mathsf{s}'}}|f_\mathsf{s}-f_\mathsf{s'}|\le C_1\|f\|_{C^l}\;\mathsf{s}^l
{\sup_{\substack{\a\in\nn\\ |\a|_1\le l'}}\sup_{{\mathcal{O}_{\mathsf{s}'}}}|\dpr^\a f_\mathsf{s}-\dpr^\a f_\mathsf{s'}|\le C_1\|f\|_{C^l}\;\mathsf{s}^{l-l'}},\quad 
\|f-f_\mathsf{s}\|_{C^{l'}}\le C_1\|f\|_{C^l}\;\mathsf{s}^{l-l'},
\eeq 
for any $0<\mathsf{s}'<\mathsf{s}$ and any $0\le l'\le l$ with $l'\in\natural$. If, in addition, $f$ is periodic in some component $y_j$ or $x_j$, then so is $f_\mathsf{s}$ in that component.
\elem

\lemtwo{BernsMos}{Bernstein, Moser}
Assume that $\{f_j\}_{j\ge 0}$ is a sequence of real--analytic functions defined respectively on $\mathcal{O}_j\coloneqq \{(y,x)\in \cn\times\cn\,:\,|\Im(y,x)|<\mathsf{s}_j\}$ such that
$$
\sup_{\mathcal{O}_j}|f_j-f_{j-1}|\le \g\;\mathsf{s}_{j-1}^{l},\qquad \forall\; j\ge 1,
$$
where $l\in\real_+\setminus \integer$, $\g>0$ and $\mathsf{s}_j\coloneqq \mathsf{s}_0\x^j$, with $\mathsf{s}_0>0$ and $0<\x<1$. %$\{\mathsf{s}_j\}_{j\ge 0}$ is a decreasing sequence of positive real--number with $\sum_{j\ge 0}\mathsf{s}_j<\infty$. 
Then, $f_j$ converges uniformly on $\rn$ to a function $f\in C^l(\rn\times\tn)$. %, with
%\beq{MosBerns}
%\|f-f_0\|_{C^l}\le \g C_2 ,
%\eeq
%for some constant $C_2=C_2(d,l)$. 
Moreover, if all the $f_j$ are periodic in some component $y_i$ or $x_i$, then so is $f$ in that component.
\elem
\subsection{Iteration of the KAM step and convergence}
Let $\mathsf{K},\,\mathsf{T},\,\eta ,\,\vae,\, \s,\,\r,\, \a_*,\, \b_0 $ be as in \S\ref{settings} and \ref{SecStatem}. 
Let  $0<m<l/2-\n$, $0<\wh m<\min\{(m+1)/\n,2\}$, $l'\coloneqq \max\{(6+2l/\n)(l+\n)/(l-2\n)-2l(l-\n)/(l-2\n),\; 2l(l+3\n)/(\n(l-2\n))\}$  and for $j\ge 0$, let
\begin{align*}
&\s_0\coloneqq C_2^{-1}\s\;,\quad s_{0}\coloneqq  4\s_0\;,\quad r_0\coloneqq \a\s_0^\n /(2\mathsf{K})\;,  \quad  \l\coloneqq \log\r^{-1}\;,\\
& \x\coloneqq (C_2\eta^{1/\n}\l)^{-1},\quad\s_j\coloneqq \s_0\x^j,\quad s_{j}\coloneqq 4\s_j= 4\s_0\x^j,\quad \bar{\s}_j\coloneqq \x^{mj},\quad \k_j\coloneqq 6\s_j^{-1}\l,\\
& r_j\coloneqq r_0\x^{\n j}\;,\quad \check{r}_{j+1}\coloneqq \frac{r_0}{64d\eta}\x^{\n j}\;, \quad \tilde{r}_{j+1}\coloneqq  \frac{r_0}{2^{11}d^2\eta^2}\x^{(\n+m) j}\;,\quad \x_0\coloneqq s_0\;,\quad\x_{j+1}\coloneqq \s_j\;,\\
&\mathcal{S}_j\coloneqq \{y\in \cn\,:\,|\Im(y)|<\x_j\}\;,\quad \mathcal{O}_j\coloneqq \{(y,x)\in \cn\times(\cn/\zn)\,:\,|\Im(y,x)|<\x_j\}\;,\\
&\|\cdot\|_{\x_j}\coloneqq \sup_{\mathcal{O}_j}|\cdot|\;.
\end{align*}
First of all, we extend $\mathrm{K}$ and ${\mathrm{P}}$ to the whole phase space $\rn\times\tn$. 
\subsubsection{Extension of $K$ and $P$ to the whole space}
%In the present section, we aim to extend $\mathrm{K}$ and ${\mathrm{P}}$ to $\rn\times\tn$. 

\noi
First of all, there exist\footnote{see for instance \cite[Lemma~2.2.1]{koudjinan2019quantitative}} $\mathsf{C}_0=\mathsf{C}_0(d,l)>0$ and a {\it Cut--off} $\chi\in C(\cn)\cap C^\infty(\rn)$ with $0\le \chi\le 1$, $\supp\chi\subset D_{\a_*}(\mathscr D')$, $\chi\equiv 1$ on $D_{\a_*/2}(\mathscr D')$ and for any $k\in\nn$ with $|k|_1\le l$,
$$
\|\dpr_y^k \chi\|_{\rn}\le \mathsf{C}_0\;\a_*^{-|k|_1}\;.
$$
%where $\mathscr D' \coloneqq \left\{y_0\in\mathscr D:  B_{2\a}(y)\subseteq\mathscr D\right\}$.
 By the F\`aa Di Bruno's Formula\cite{CS96Mult}, there exists $\mathsf{C}_1=\mathsf{C}_1(d,l)>0$ such that for any $f\in C^l(\rn\times\tn)$, we have \beq{FaaMultC1}
\| \chi\circ\pi_1\cdot f\|_{C^l}\le \mathsf{C}_1\;\a_*^{-l}\;\|f\|_{C^l}.
\eeq
\noi
Let $ \wh{\mathrm{K}}\in C^l(\rn)$ such that\footnote{Observe that $\mathsf{C}_1^{-1}\a_*^{l}/4<1/4$.} $\|T\|_{\mathscr D}\|\wh{\mathrm{K}}-\mathrm{K}\|_{C^l(\mathscr D)}\le \mathsf{C}_1^{-1}\a_*^{l}/4$. Thus, $\wh{\mathrm{K}}_{yy}= \mathrm{K}_{yy}(\uno_d+T(\wh{\mathrm{K}}_{yy}-\mathrm{K}_{yy}))$ is invertible on $\mathscr D$ and $\|(\wh{\mathrm{K}}_{yy})^{-1}\|_{\mathscr D}\le 2\|T\|_{\mathscr D}$.  Then, $K\coloneqq \wh {\mathrm{K}}+\chi\cdot ({\mathrm{K}}-\wh {\mathrm{K}})\in C^l(\rn\times\tn)$, $K\equiv {\mathrm{K}}$ on $D_{\a_*/2}(\mathscr D')$ and %using the F\`aa Di Bruno's Formula, we get
$$
\| K\|_{C^l}\leby{FaaMultC1} \| {\mathrm{K}}\|_{C^l}+\mathsf{C}_1\;\a_*^{-l}\;\|\wh {\mathrm{K}}-{\mathrm{K}}\|_{C^l}\le \|{\mathrm{K}}\|_{C^l}+\|T\|_{C^l}^{-1}/4<  2\|{\mathrm{K}}\|_{C^l}
$$
and
$$
\| (\wh {\mathrm{K}}_{yy})^{-1}\dpr^2_y(\chi\cdot ({\mathrm{K}}-\wh {\mathrm{K}}))\|_{\mathscr D}\le \|(\wh {\mathrm{K}}_{yy})^{-1}\|_{\mathscr D}\cdot\mathsf{C}_1\;\a_*^{-l}\;\|\wh {\mathrm{K}}-{\mathrm{K}}\|_{C^l}\le 1/2.
$$
Therefore, $K_{yy}$ is in particular invertible and $\|(K_{yy})^{-1}\|_{\mathscr D}\le 2 \|(\wh {\mathrm{K}}_{yy})^{-1}\|_{\mathscr D}\le 4\|T\|_{\mathscr D}$.

\noi
Similarly, one extends ${\mathrm{P}}$ to a function $P\in C^l(\rn\times\tn)$ such that $K\equiv {\mathrm{K}}$ on $D_{\a_*/2}(\mathscr D')$ and $\| P\|_{C^l}\le  2\|{\mathrm{P}}\|_{C^l}$. Now, letting $H\coloneqq K+  P$, we have $H_{|D_{\a_*/2}(\mathscr D')}={\mathrm{H}}$. Hence, it does not make any difference for us replacing ${\mathrm{H}}$ by $H$ since  the invariant tori of $H$ we shall construct live precisely in $D_{\a_*/2}(\mathscr D')$ as $r_0\ltby{SmaLConD} \a_*/2$.

\ \\
\noi
Let $\mathcal{K}_j$ (resp. $\mathcal{P}_j$) be the real--analytic approximation $K_{\x_j}$ (resp. $P_{\x_j}$) of $K$ (resp. $P$) defined on $\mathcal{O}_j$ given by Lemma~\ref{JacMosZen}.  Then, the following holds.
%%Then, set
%%\begin{align*}
%%&C_8\coloneqq \max\{a^{l-l'}C_1,\}\;,\\
%%&C_7\coloneqq \max\{2\sqrt{2} C_1,\}\;,\\
%%\end{align*}

\subsubsection{Iteration of the KAM Step}
There exist  constants $C_j=C_j(d,\t,l)>1$ ($j=1,\cdots,6$) such that the following holds.
\lem{IteKAM} %Let $2\le l'< l-1$ and
Set 
$\mathscr D_0\coloneqq
%(\dpr_{y}\mathcal{K}_0)^{-1}(\dpr_y K(\mathscr D_{\a}))\overset{def}{=}
\{y\in\rn\ :\ \dpr_{y}\mathcal{K}_0(y)\in \dpr_y K(\mathscr D_{\a})\}$.
Assume that
\beq{SmaLConD}
\a\le C_3^{-1}  \mathsf{K}\qquad \mbox{and}\qquad C_3\;\vae\;\mathsf{K}^{\frac{l+2\n}{l-2\n}}\eta^{l'}\a^{\frac{-2l}{l-2\n}}\le 1\;.
\eeq
Then, the following assertions $(\mathscr P_j)$, $j\ge 1$, hold. There exist a sequence  of sets $\mathscr D_j$, a sequence of diffeomorphisms $G_j\colon D_{\tilde{r}_j}(\mathscr D_{j-1})\to G_j(D_{\tilde{r}_j}(\mathscr D_{j-1}))$, a sequence of real--analytic symplectic transformations%\footnote{Observe that $s_j+\s_{j-1}/3=(3\x+1)\s_{j-1}/3\ltby{SmaLConD} 2\s_{j-1}/3$ and $2r_j+  r_{j-1}\s_{j-1}/3\ltby{SmaLConD} \s_0^{\n}\x^{j}+ \s_{j-1}/3=\s_0^\t\s_j+\s_{j-1}/3< \s_{j-1}$, which combined with \equ{phiokExt0Ext01} imply \equ{phijBisv2v501}.}
\beq{phijBisv2v501}
\phi_{j}=(v_j,u_j):D_{r_{j},s_{j}}(\mathscr D_{j})\to D_{\s_{j-1},\s_{j-1}}(\mathscr D_{j-1})\;,
\eeq
such that, setting $\mathcal{H}_{j-1}\coloneqq \mathcal{K}_{j-1} +  \mathcal{P}_{j-1}$, we have 
\begin{align}
&G_j(\mathscr D_{j-1})=\mathscr D_j\subset \mathscr D_{r_j} \;, \qquad G_{j}=(\dpr_y K_{j})^{-1}\circ \dpr_y K_{j-1}\;,\label{FinGjj}\\
%&\dpr_{y}K_{j}\circ G^{j}=\dpr_{y}\mathcal{K}_{j-1} \;,\\
%&\phi_{j+1}:D^d_{r_{j+1},s_{j+1}}(\mathscr D_{j+1})\to D^d_{r_{j},s_{j}}(\mathscr D_{j})\quad\qquad\quad\mbox{is real--analytic},\label{phijExt}\\
&\det\dpr_y^2 K_j(y)\neq 0,\quad T_j(y)\coloneqq \dpr_y^2 K_j(y)^{-1}\;,\hspace{1.55cm}\forall\; y\in \mathscr D_j\;,\label{FinDEtkJ}\\
&H_{j}\coloneqq \mathcal{H}_{j-1}\circ\phi^{j}\eqqcolon K_{j} +  P_{j}\hspace{3.5cm}\  \mbox{on}\quad D_{r_{j},s_{j}}(\mathscr D_{j})\;,\label{KAMToriH}
&  %G_1=(\dpr_y K_1)^{-1}\circ \dpr_y \mathcal{K}_0 \quad\mbox{and}\quad G_{j+1}=(\dpr_y K_{j+1})^{-1}\circ \dpr_y K_{j}\;,
\end{align}
	%\beq{HjExt}
	%H_j:=H_{j-1}\circ\phi_j=: K_j + \vae^{2^j} P_j\ ,
	%\eeq
%converge uniformly, 
where $\phi^j\coloneqq \phi_1\circ \phi_2\circ \cdots\circ \phi_j$ and $K_0\coloneqq \mathcal{K}_0$.% and $G^j\coloneqq G_{j}\circ G_{j-1}\circ\cdots\circ G_2\circ G_1$. %More precisely, we have the following:

%\textcolor{red}
{%Finally, the following estimates hold for any $i\ge 1$:%\footnote{Observe that \equ{estGidevidv2} follows \equ{estGiidv2} using Cauchy's estimate.}
%\beqa{estfin1Ext}
%&& 
Moreover,%\footnote{In principle, instead of $\phi_j-\id$, one should write $\phi_j-(\pi_1,\pi_2)$. But, for the sake of simplicity we shall use this abuse of notation. Hence, $v_j-\id$ should be understood $v_j-\pi_1$ , \etc.}
\begin{align}
&\|G_j-\id\|_{\tilde r_{j},\mathscr D_{j-1}}\le 
\tilde{r}_j\;\x^{2\n} {\x^{m(j-1)}}\;,\qquad \|\dpr_y G_j-\uno_d\|_{\tilde r_{j},\mathscr D_{j-1}}\le 
\x^{2\n} {\x^{2m(j-1)}},\label{estGiidv201}\\
%&\|\dpr_z G_{j}-\uno_d\|_{\tilde{r}_{j}/2,\mathscr D_{j-1}}\le 64d\eta_0\;\s_{j-1}^{\n+d}\;|\vae|^{2^{j-1}}\mathsf L_{j-1} \;,\label{estGidevidv2}\\
&\|\dpr_y^2 K_j\|_{ r_{j},\mathscr D_{j}}< 
2\mathsf{K}\;,\qquad\quad \|T_j\|_{ \mathscr D_{j}}< 2\mathsf{T}\;,\qquad\quad T_{j}\coloneqq (\dpr_y^2 K_j)^{-1}\;,\label{estGi3001}\\%\label{estGi3101}\\
&\|P_j\|_{r_j,s_j,\mathscr D_j}\le  C_1\;\mathsf{K} \; \x_{j-1}^l\ ,
\label{estfin2Ext01v501}\\
&\max\big\{\|\mathsf{W}_j(\phi_j-\id)\|_{2r_j,s_j,\mathscr D_j}\;,\; \|\pi_2\dpr_x(\phi_j-\id)\|_{2r_j,s_j,\mathscr D_j}\big\}
\le  \x^{2\n}\x^{m(j-1)}\label{estfin2Ext03v501} ,
\end{align}
} 
where $\mathsf{W}_j\coloneqq \diag(r_{j-1}^{-1}\uno_d,\; \s_{j-1}^{-1}\uno_d)$.
\elem
\rem{DomainFitness}
Observe that 
\begin{align}
&s_j+\s_{j-1}/3=(12\x+1)\s_{j-1}/3\ltby{SmaLConD} 2\s_{j-1}/3 <s_{j-1}/2\;,\label{FitnessEq1}\\
& 2r_j+  r_{j-1}\s_{j-1}/3\ltby{SmaLConD}r_{j-1}/4+r_{j-1}/6<r_{j-1}/2\;,\label{FitnessEq2}\\
& 2r_j+  r_{j-1}\s_{j-1}/3\ltby{SmaLConD} \s_0^{\n}\x^{j}+ \s_{j-1}/3=\s_0^\t\s_j+\s_{j-1}/3< \s_{j-1}\;,\label{FitnessEq3}
\end{align} which combined with \equ{phiokExt0Ext01} imply
\beq{VIPRel}
\phi_{j}(D_{r_{j},s_{j}}(\mathscr D_{j}))\subset D_{\s_{j-1},\s_{j-1}}(\mathscr D_{j-1})\bigcap D_{r_{j-1}/2, s_{j-1}/2}(\mathscr D_{j-1})\;,
\eeq
and, in particular, \equ{phijBisv2v501}. %Notice also that, \equ{EqRJPl1} will allow us to apply Cauchy's estimate to \equ{convEstExt01} and bound the differential $D\phi^j$ also on $D_{r_{j},s_{j}}(\mathscr D_{j})$.
Also,
\beq{EqRJPl1}
2r_{j+1}\ltby{SmaLConD} \su 4 \dst\min\left\{\frac{\a}{2d(2\mathsf{K})\k_j^{\n}}\,,\, \check{r} _{j+1}\right\}\;,
\eeq
which, together with the definitions of the sequences of the various parameters, implies that \equ{DefNArn2v501} is fulfilled for any $j\ge 1$.
\erem
%\subsubsection{Proof Iteration}
\proof\\
{\bf Step 1: }We check $(\mathscr P_1)$. We claim that we can apply Lemma~\ref{lem:1Extv5Simpl01} to $\mathcal{H}_0$. Indeed, we have%\footnote{Observe that $r_0= s_0=\x_0.$}
\beq{HessApprK0}
\|\mathcal{P}_0\|_{r_0}\le \|\mathcal{P}_0\|_{\x_0}\leby{AnalAppro} C_1\|P\|_{C^0}\le 2C_1\vae.
\eeq
 From
$$
\dpr^2_y\mathcal{K}_0=\dpr^2_y K(\uno_d+T\dpr^2_y(\mathcal{K}_0-K))
$$
and
$$
\|T\dpr^2_y(\mathcal{K}_0-K)\|_{r_0,\mathscr D_0}\le\sup_{\mathcal{S}_0}\|T\dpr^2_y(\mathcal{K}_0-K)\|
%\le \mathsf{T}\|\dpr^2_y(\mathcal{K}_0-K)\|_{C^{l'}}
\leby{AnalAppro}C_1\eta\;s_0^{l-2}\le \su2\;,
$$
it follows that $\dpr^2_y\mathcal{K}_0$ is invertible, $\|\dpr^2_y\mathcal{K}_0\|_{r_0,\mathscr D_0}\le 2\mathsf{K}$
 and 
\beq{ApprK0}
\|(\dpr^2_y\mathcal{K}_0)^{-1}-\mathsf{T}\|_{\mathscr D_0}\le 2\mathsf{T}C_1\eta\;s_0^{l-2}<\mathsf{T}\;,\qquad \|(\dpr^2_y\mathcal{K}_0)^{-1}\|_{\mathscr D_0}< 2\mathsf{T}.
\eeq
Hence, $\mathcal{H}_0$ satisfies the assumptions in \equ{RecHypArnExtv501} with\footnote{``$a\leadsto b$'' stands for ``$a$ replaced by $b$''.} $\vae\leadsto C_1\vae$, $r\leadsto r_0$, $s\leadsto s_0$, $\s\leadsto \s_0$, $\mathsf{K}\leadsto 2\mathsf{K}$. Now, observe that 
$$
\r=\frac{2C_1\vae}{\a r_0\s_0^\n}\leby{SmaLConD} \su4 \qquad \mbox{and}\qquad r_0=\frac{\a}{2\mathsf{K}}\s_0^\n,
$$
and, therefore, \equ{DefNArnExt1v501} is verified. Moreover, by the definitions and \equ{EqRJPl1}, \equ{DefNArn2v501} holds trivially with $\mathsf{T}\leadsto 2\mathsf{T}$, $\eta\leadsto 4\mathsf{T}\mathsf{K}$ $\k\leadsto\k_0$, $\check{r}\leadsto\check{r}_1$, $\bar{r}\leadsto4r_1$, $\tilde{r}\leadsto\tilde{r}_1$, $\bar{s}\leadsto s_0-2\s_0/3$, $s'\leadsto s_1<s_0-\s_0$, $\bar{\s}\leadsto \bar{\s}_0$, $\mathsf{L}\leadsto 32 C_0C_1 \mathsf{K}\mathsf{T}^2\vae/(r_0\tilde{r}_1)\leby{SmaLConD} \bar{\s}_0/3$.
 Consequently, we can apply Lemma~\ref{lem:1Extv5Simpl01} to $\mathcal{H}_0$ to obtain  the change of coordinates $\phi^1=\phi_1$. In particular,   \equ{phiokExt0Ext01} yields ${\equ{phijBisv2v501}}_{j=1}$, \equ{HPhiH'01} yields ${\equ{FinGjj}}_{j=1}\div{\equ{KAMToriH}}_{j=1}$ and, \equ{convEstExt01} yields ${\equ{estGiidv201}}_{j=1}\div{\equ{estfin2Ext03v501}}_{j=1}$.
Therefore, $(\mathscr P_1)$ is proven. \\
% with $r=r_0$ $\mathsf{L}=C_0C_1\frac{\eta\mathsf{T}\vae}{r_0\tilde{r}_1}$.\\

%\subsubsection{Step 2}
\noi
{\bf Step 2: }We assume $(\mathscr P_j)$  holds for some $j\ge 1$ and check $(\mathscr P_{j+1})$. Write
$$
\mathcal{H}_{j}\coloneqq \mathcal{K}_{j} +  \mathcal{P}_{j}=\mathcal{H}_{j-1}+ (\mathcal{K}_{j}-\mathcal{K}_{j-1}) +  (\mathcal{P}_{j}-\mathcal{P}_{j-1}).
$$
%%Observe first that
%%\beq{2rjSgmjMoins1}
%%3r_j=3r_0\x^{\b j}\le 3r_0\frac{\s_{j-1}}{\s_0}\x^{\b}< 3\mathsf{K}\s_0^{\b-1}\s_{j-1}\le \s_{j-1}=\x_j.
%%\eeq
By the inductive assumption and \equ{VIPRel}, we have%\marginnote{$r_{j-1}\s_{j-1}/4\le\x_{j}$.}
\begin{align*}
\mathcal{H}_{j}\circ\phi^j &=\mathcal{H}_{j-1}\circ\phi^j+ (\mathcal{K}_{j}-\mathcal{K}_{j-1})\circ v^j +  (\mathcal{P}_{j}-\mathcal{P}_{j-1})\circ\phi^j\\
   &= K_{j} +   P_{j}+ (\mathcal{K}_{j}-\mathcal{K}_{j-1})\circ v_j +  (\mathcal{P}_{j}-\mathcal{P}_{j-1})\circ\phi^j\\
   &=\mathcal{K}^j+\mathcal{P}^j\hspace{1.5cm}\  \mbox{on}\quad D_{r_{j},s_{j}}(\mathscr D_{j})\;,
\end{align*}
where $\mathcal{K}^j\coloneqq K_{j}$ and $\mathcal{P}^j\coloneqq  P_{j}+(\mathcal{K}_{j}-\mathcal{K}_{j-1})\circ\phi^j+   (\mathcal{P}_{j}-\mathcal{P}_{j-1})\circ\phi^j$, with
\beq{EqmAtCPjTj}
\|\dpr_y^2\mathcal{K}^j\|_{r_j,\mathscr D_j}< 2\mathsf{K},\qquad \|(\dpr_y^2\mathcal{K}^j)^{-1}\|_{\mathscr D_j}< 2\mathsf{T},
\eeq
by the inductive assumption, provided $\phi^j$ maps $D_{r_{j},s_{j}}(\mathscr D_{j})$ into $\mathcal{O}_j=\{(y,x)\in \cn\times\cn\,:\,|\Im(y,x)|<\x_j\}$ 
\ie
\beq{TTiMp}
\sup_{ D_{r_{j},s_{j}}(\mathscr D_{j})}|\Im\phi^j|\le \frac{\x_j}{2} ,
\eeq 
%Indeed, letting $\dpr^2_y\mathcal{K}^j= \dpr^2_y K_{j}(\uno_d+T_j\dpr^2_y(\mathcal{K}_{j}-\mathcal{K}_{j-1}))\eqqcolon \dpr^2_y K_{j}(\uno_d+A)$, we have, using Cauchy's estimate, %under \equ{TTiMp}
%\begin{align*}
%\sup_{ D_{r_{j},s_{j}}(\mathscr D_{j})}\|A\|&\le 8\mathsf{T}\|\dpr_y^2(\mathcal{K}_{j}-\mathcal{K}_{j-1})\|_{\x_j,\mathscr D_j}\leby{AnalAppro} 8C_1\eta \x_j^{l-2}\\
%	&\le   C_3\eta\s_0^{l-2}\x^{( l-2)j}\le  C_3\eta\s_0^{l-2}\le\su2, %\marginnote{$\s_0\le (C_3\eta)^{-1/{(l-2)}}$}
%\end{align*}
%and \equ{EqmAtCPjTj} follows at once. Now, we prove \equ{TTiMp}. 
which we now prove.
 Observe that, for any $1\le i\le j$,\footnote{$\|\mathsf{W}_i(D\phi_i-\uno_{2d})\mathsf{W}_i^{-1}\|_{r_i,s_i,\mathscr D_i}=\max\{\|\dpr_y v_j-\Id\|_{r_i,s_i,\mathscr D_i},\frac{\s_{j-1}}{r_{j-1}}\|\dpr_x v_j\|_{r_i,s_i,\mathscr D_i},\frac{r_{j-1}}{\s_{j-1}}\|\dpr_y u_j\|_{r_i,s_i,\mathscr D_i},\|\dpr_x u_j-\Id\|_{r_i,s_i,\mathscr D_i}\}$.}
\beq{DevphiJ}
\|\mathsf{W}_i\mathsf{W}_{i+1}^{-1}\|%=\max\{\x^\n,\x\}
=\x
\quad \mbox{and} \quad\|\mathsf{W}_i(D\phi_i-\uno_{2d})\mathsf{W}_i^{-1}\|_{r_i,s_i,\mathscr D_i}\leby{estfin2Ext03v501} 2\x^{\n}\x^{m(i-1)}.
\eeq
Thus, writing $\mathsf{W}_1D\phi^j\mathsf{W}_j^{-1}=(\mathsf{W}_1D\phi_1\mathsf{W}_1^{-1})(\mathsf{W}_1\mathsf{W}_{2}^{-1})\cdots(\mathsf{W}_jD\phi_j\mathsf{W}_j^{-1})$, we then get from \equ{DevphiJ}:
\begin{align}
\|\mathsf{W}_1D\phi^j\mathsf{W}_j^{-1}\|_{r_j,s_j,\mathscr D_j}&\le \x^{j-1}\prod_{i=1}^j(1+2\x^{\n}\x^{m(i-1)})\nonumber\\
	&\le \x^{j-1}\exp(4\x^{\n})\leby{SmaLConD} 2\x^{j-1}.\label{hihih}
\end{align}
Now, writing $(y,x)= X+\mathbf{i}Y$, where $X=X(y,x)\coloneqq \Re (y,x)\in \rn\times\rn$ and $Y=Y(y,x)\coloneqq \Im (y,x)\in \rn\times\rn$, %are both elements of $\rn\times\rn$ 
we obtain
$$
\phi^j(y,x)= \phi^j(X)+\mathbf{i}\mathscr I(y,x),
$$
where
$$
\mathscr I(y,x)\coloneqq \mathsf{W}_1^{-1}\cdot\int_0^1\mathsf{W}_1 D \phi^j(X+\mathbf{i}tY)\mathsf{W}_j^{-1}dt\cdot \mathsf{W}_jY,
$$
so that, as $\phi^j$ is real on reals, we have $\phi^j(X)\in\rn$ and, therefore, $\Im (\phi^j(y,x))= \Re(\mathscr I(y,x)).$ 
We have, for any $y_j\in \mathscr D_{j}\subset \rn$ and $y\in D_{r_{j}}(y_{j})$,%\footnote{Recall that $\mathscr D_0\subset\rn$.}
%\begin{align*}
%\frac{|\Im(y)|}{r_{j-1}}&\le \su{r_{j-1}}\left(|\Im(y_0)|+|\Im(y-y_j)|+\sum_{i=0}^{j-1}|\Im(y_{i+1}-y_i)|\right)\\
%	&\le \su{r_{j-1}}\left(r_j+\sum_{i=0}^{j-1}\tilde{r}_{i+1}\right)
%\end{align*}
\beq{eqDomFit1}
\frac{|\Im(y)|}{r_{j-1}}\le \su{r_{j-1}}\left(|\Im(y-y_j)|+|\Im(y_j)|\right)\le \su{r_{j-1}}|y-y_j|\le \su{r_{j-1}}r_j=\x^\n\;,
\eeq
so that
\beq{eqDomFit2}
\sup_{ D_{r_{j},s_{j}}(\mathscr D_{j})}|\mathsf{W}_j Y|=\sup_{(y,x)\in D_{r_{j},s_{j}}(\mathscr D_{j})}\max\left\{\frac{|\Im(y)|}{r_{j-1}}\;,\;\frac{|\Im(x)|}{\s_{j-1}}\right\}\leby{eqDomFit1} \max\left\{\x^\n\;,\;4\x\right\}=4\x.
\eeq
Moreover, since $\|\mathsf{W}_1^{-1}\|\eqby{SmaLConD}\s_0$, we have, for any $(y,x)\in D_{r_{j},s_{j}}(\mathscr D_{j})$,
\begin{align*}
|\Im (\phi^j(y,x))|&=|\Re(\mathscr I(y,x))|\le \|\mathsf{W}_1^{-1}\| \left\|\int_0^1\mathsf{W}_1 D \phi^j(X+\mathbf{i}tY)\mathsf{W}_j^{-1}dt\right\| |\mathsf{W}_jY|\\
	&\stackrel{\equ{hihih}+\equ{eqDomFit2}}{\le}\s_0\cdot 2\x^{j-1}\cdot 4\x=8\s_{j}<\x_j,
\end{align*}
which completes the proof of \equ{TTiMp}. Hence,\footnote{It may seem artificial treating $(\mathcal{K}_{j}-\mathcal{K}_{j-1})\circ \pi_1\circ\phi^j$ as a reminder, as the latter is of order $1$ compared to $(\mathcal{P}_{j}-\mathcal{P}_{j-1})\circ\phi^j$. Indeed, to get \equ{calPjxijmoin1}, we have bounded $\vae$ by $\mathsf{K}$, and the latter is of order $1$. The point is that the norm of $\phi^j-\id$ is in fact large, it is basically of order $\|\phi^1-\id\|_{r_1,s_1,\mathscr D_1}$. Hence, extracting the integrable part of $(\mathcal{K}_{j}-\mathcal{K}_{j-1})\circ \pi_1\circ\phi^j$ does not yield significant improvement. %On the other hand, this may be the origin of the  relatively larger measure estimate of the Kolmogorov's set we obtain at the end, but we were unable to do better.
}
\begin{align}
\|\mathcal{P}^j\|_{r_j,s_j,\mathscr D_j}&\le \|  P_{j}\|_{r_j,s_j,\mathscr D_j}+\|(\mathcal{K}_{j}-\mathcal{K}_{j-1})\circ\phi^j \|_{r_j,s_j,\mathscr D_j}+   \|(\mathcal{P}_{j}-\mathcal{P}_{j-1})\circ\phi^j\|_{r_j,s_j,\mathscr D_j}\nonumber\\
   &\stackrel{\equ{phijBisv2v501}+\equ{estfin2Ext01v501}+\equ{TTiMp}}{\le} C_1 \mathsf{K}{\x_{j-1}^{l}}+\|\mathcal{K}_{j}-\mathcal{K}_{j-1} \|_{\x_j}+  \|\mathcal{P}_{j}-\mathcal{P}_{j-1}\|_{\x_j}\nonumber\\ 
   &\leby{AnalAppro} C_1 \mathsf{K}{\x_{j-1}^{l}}+C_1 \mathsf{K}{\x_{j-1}^{l}}+C_1\;\vae\;\x_{j-1}^l \nonumber\\
   &\leby{SmaLConD} 3C_1 \mathsf{K}{\x_{j-1}^{l}}\;.\label{calPjxijmoin1}
\end{align}
Thus, thanks to \equ{EqmAtCPjTj}, $\mathcal{H}_{j}\circ\phi^j=\mathcal{K}^j+\mathcal{P}^j$ satisfies the assumptions in \equ{RecHypArnExtv501} with $\vae\leadsto \|\mathcal{P}^j\|_{r_j,s_j,\mathscr D_j}$, $r\leadsto r_j$, $s\leadsto s_j$, $\s\leadsto \s_j$, $\mathsf{K}\leadsto 2\mathsf{K}$ as 
$$
\dpr^2_y\mathcal{K}^j(\mathscr D_j)\overset{def}=\dpr^2_y K_j(G_j(\mathscr D_{j-1}))=\dpr^2_y K_{j-1}(\mathscr D_{j-1})=\cdots=\dpr^2_y K_{0}(\mathscr D_{0})\subset \D^\a_\t.
$$
%Next, observe that 
%$$
%\r=\frac{2C_1\vae}{\a r_0\s_0^\n}\leby{SmaLConD} \su4 \qquad \mbox{and}\qquad r_0=\frac{\a}{2\mathsf{K}}\s_0^\n,
%$$
%and, therefore, \equ{DefNArnExt1v501} is verified. Moreover, by the definitions and \equ{EqRJPl1}, \equ{DefNArn2v501} holds trivially with

\noi
Hence, in order to apply Lemma~\ref{lem:1Extv5Simpl01} to $\mathcal{H}_{j}\circ\phi^j=\mathcal{K}^j+\mathcal{P}^j$, we need only to check \equ{DefNArnExt1v501}, \equ{DefNArn2v501} and \equ{cond1ExtExtv501}. But, by the definitions and \equ{EqRJPl1}, \equ{DefNArn2v501} holds trivially with 
 $\mathsf{T}\leadsto 2\mathsf{T}$, $\eta\leadsto 4\mathsf{T}\mathsf{K}$ $\k\leadsto\k_j$, $\check{r}\leadsto\check{r}_{j+1}$, $\bar{r}\leadsto 4r_{j+1}$, $\bar{\s}\leadsto \bar{\s}_j$, $\tilde{r}\leadsto\tilde{r}_{j+1}$, $\bar{s}\leadsto s_j-2\s_j/3$, $s'\leadsto s_{j+1}<s_{j}-\s_{j}$, $\mathsf{L}\leadsto 32 C_0 \mathsf{K}\mathsf{T}^2\|\mathcal{P}^j\|_{r_j,s_j,\mathscr D_j}/(r_j\tilde{r}_{j+1})$. %$\leby{SmaLConD} \bar{\s}_0/3$. 
 Moreover, we have, 
$$
r_j= r_0\x^{\n j}\le \frac{\a}{2\mathsf{K}}\s_j^\n\leby{EqmAtCPjTj}
{\a\s_j^\n}/{\|\dpr_y^2\mathcal{K}^j\|_{r_j,s_j,\mathscr D_j}}\;,
$$
%%so that
%%$$
%%\s_1^{-\n}\frac{\|\mathcal{P}^1\|_{r_1,s_1,\mathscr D_1}}{{\a}r_1}\le 3C_1\mathsf{K}\s_0^{-\n}\x^{-(\n+\n)}s_0^l(\a r_0)^{-1}
%%$$
%%and for $j\ge 2$,
%%$$
%%\s_j^{-\n}\frac{\|\mathcal{P}^j\|_{r_j,s_j,\mathscr D_j}}{{\a}r_j}\le 3C_1\mathsf{K}\s_j^{-\n}\;.
%%$$
%
%
\begin{align*}
\s_1^{-\n}\frac{\|\mathcal{P}^1\|_{r_1,s_1,\mathscr D_1}}{{\a}r_1}\r^{-1}&\leby{calPjxijmoin1} C_2\s_0^{l}\frac{\mathsf{K}}{\vae\x^{2\n}}\leby{SmaLConD} 1, 
\end{align*}
\begin{align*}
3C_0\frac{\eta\mathsf{T}\|\mathcal{P}^1\|_{r_1,s_1,\mathscr D_1}}{r_1\tilde{r}_{2}\bar{\s}_1 }&\leby{calPjxijmoin1}  C_2\s_0^{l-2\n}\frac{\eta^{6+m}\mathsf{K}^2}{\a^2}\l^{2(\n+m)}\leby{SmaLConD} \x^{2\n}, 
\end{align*}
and, for $j\ge 2$,
\begin{align}
\s_j^{-\n}\frac{\|\mathcal{P}^j\|_{r_j,s_j,\mathscr D_j}}{{\a}r_j}\r^{-1}&\le C_2\s_0^{l}\frac{\mathsf{K}}{\vae}\x^{(l-2\n)j-2l}\le C_2\s_0^{l}\frac{\mathsf{K}}{\vae}\x^{-2l} \leby{SmaLConD} 1\;,\label{ConvDfiJ104g}
\end{align}
\begin{align*}
3C_0\frac{\eta\mathsf{T}\|\mathcal{P}^j\|_{r_j,s_j,\mathscr D_j}}{r_j\tilde{r}_{j+1}\bar{\s}_j}&\leby{calPjxijmoin1} C_2 \s_0^{l-2\n}\frac{\eta^{4+2l/\n}\;\mathsf{K}^2}{\a^2}\l^{2l} \leby{SmaLConD} \x^{2\n},
\end{align*}
which concludes the verification of \equ{DefNArnExt1v501} and \equ{cond1ExtExtv501}. Therefore, Lemma~\ref{lem:1Extv5Simpl01} applies to $\mathcal{H}_j$ and yields the desired symplectic change of coordinates $\phi_{j+1}$. In particular,   \equ{phiokExt0Ext01} yields ${\equ{phijBisv2v501}}_{j+1}$, \equ{HPhiH'01} yields ${\equ{FinGjj}}_{j+1}\div{\equ{KAMToriH}}_{j+1}$ and, \equ{convEstExt01} yields ${\equ{estGiidv201}}_{j+1}\div %, \equ{estGi3001}_{j+1}$ and $
{\equ{estfin2Ext03v501}}_{j+1}$, as 
$$
C_1\r\cdot 3C_1\mathsf{K}\x_{j-1}^l=  C_1\mathsf{K}\x_{j}^l\cdot C_1\r(\eta^{1/\n}\log\r^{-1})^{l}\le C_1\mathsf{K}\x_{j}^l\cdot C_1\eta^{l/\n}\r^{1/2}\leby{SmaLConD}C_1\mathsf{K}\x_{j}^l.
$$
This ends the proof of  $(\mathscr P_{j+1})$, and, consequently, of the  Lemma.  %, with $P_{j+1}\coloneqq \mathcal{P}_j'$
\qed
\subsubsection{Convergence of the procedure}
Now, we are in position to prove the convergence of the KAM scheme.
\lem{ConvShc}
Under the assumptions and notation in Lemma~\ref{IteKAM}, the following holds. 
\begin{itemize}
\item[$(i)$] the sequence $G^{j}\coloneqq G_{j}\circ G_{j-1}\circ\cdots\circ G_2\circ G_1$ converges uniformly on $\mathscr D_0$ %{\d,\a}$ 
to a lipeomorphism $G_*\colon \mathscr D_{0}\to \mathscr D_*\coloneqq G_*(\mathscr D_{0})\subset\mathscr D$ and {$G_*\in C^1_W(\mathscr D_{0})$\;. } 
\item[$(ii)$]  $ P_j$ converges uniformly to $0$ on 
$\mathscr D_*\times\dst\torus^d_{s_*}$ in the $C^2_W$ topology\;;
\item[$(iii)$] $\phi^j$ converges uniformly on 
$\mathscr D_*\times\tn$ to a symplectic transformation 
$$
\phi_*\colon \mathscr D_*\times\tn\overset{into}{\longrightarrow} \mathscr D\times\tn;
$$
 with {$\phi_*\in C^{\wt m}_W(\mathscr D_*\times\tn)$ and $\phi_*(y_*,\cdot)\in C^{\wt m \n}(\tn)$, for any given $y_*\in \mathscr D_*$. %and $\phi_*(y,\cdot)\colon \torus^d_{s_*}\ni x\mapsto \phi_*(y,x)$ holomorphic, for any $y\in\mathscr D_*$
%% \quad and\quad
%%  $
%%\sup_{\mathscr D_*}|\mathsf{W}_1(\phi_*-\id)|\le 8\x^{(2\wt m-1)\n}\;.
%%$
}
\item[$(iv)$]  $K_j$ converges uniformly on 
$\mathscr D_*$ to a function {$K_*\in C^{2+\wt m}_W(\mathscr D_*)$}, with
\begin{align}
&\dpr_{y_*}K_*\circ G_*=\dpr_{y}\mathcal{K}_0 \quad \qquad\qquad\quad\quad\mbox{on} \quad {\mathscr D_0}
%\overset{????}{=}(\dpr_{y}\mathcal{K}_0)^{-1}(\dpr_y K(\mathscr D_{\a}))}
\;,\label{EqkStGSt}\\
%\mbox{so that} \quad \dpr_{y_*}K_*(\mathscr D_*)=\dpr_{y}K(\mathscr K_0)\;,\\
 &H\circ\phi_*(y_*,x)= K_*(y_*)\;,\quad\qquad\quad\  \forall(y_*,x)\in\mathscr D_*\times\tn\;.\label{EqHStfiStKst}
\end{align}
%\item[$(v)$]
\end{itemize}
\elem
%\subsubsection{Proof Conv Lemma}
\proof
The proof is essentially the same as for \cite[Lemma~6.3.3, page.~167]{koudjinan2019quantitative} , which, in turn, is based on \cite[Lemma~E.2, page.~207]{koudjinan2019quantitative}. For the reader's convenience, we give the proof for $\phi^j$; the proofs for $G^j$ and $P_j$ are similar. 
\begin{comment}
First of all, observe that, for any $j\ge 1$,
\beq{DevphiJ}
\|\mathsf{W}_j\mathsf{W}_{j+1}^{-1}\|%=\max\{\x^\n,\x\}
=\x
\quad \mbox{and} \quad\|\mathsf{W}_j(D\phi_j-\uno_{2d})\mathsf{W}_j^{-1}\|_{r_j/2,s_j/2,\mathscr D_j}\leby{estfin2Ext03v501} 2\x^{2(\wt m-1)\n}\x^{m(j-1)}.
\eeq
Thus, observing $\mathsf{W}_1D\phi^j\mathsf{W}_j^{-1}=(\mathsf{W}_1D\phi_1\mathsf{W}_1^{-1})(\mathsf{W}_1\mathsf{W}_{2}^{-1})\cdots(\mathsf{W}_jD\phi_j\mathsf{W}_j^{-1})$, we then get from \equ{DevphiJ}:
\begin{align}
\|\mathsf{W}_1D\phi^j\mathsf{W}_j^{-1}\|_{r_j/2,s_j/2,\mathscr D_j}&\le \x^{j-1}\prod_{i=1}^j(1+2\x^{2(\wt m-1)\n}\x^{m(i-1)})\nonumber\\
	&\le \x^{j-1}\exp(4\x^{2(\wt m-1)\n})\leby{SmaLConD} 2\x^{j-1},\label{hihih}
\end{align}
so that, 
\end{comment}
%Set
%\begin{align*}
%&\tilde{v}_j\coloneqq \dpr_x\pi_1(\phi_j-\id),\qquad v_j\coloneqq \id+\tilde{v}_j,\qquad V^j\coloneqq\dpr_x\pi_1\phi^j=(\uno_d+\tilde{u}_1)\circ\cdots\circ(\uno_d+\tilde{u}_j),\\
%&\tilde{u}_j\coloneqq \dpr_x\pi_2(\phi_j-\id),\qquad u_j\coloneqq \id+\tilde{u}_j,\qquad U^j\coloneqq\dpr_x\pi_2\phi^j=(\uno_d+\tilde{u}_1)\circ\cdots\circ(\uno_d+\tilde{u}_j).
%\end{align*}
Writing $\phi^{j}-\phi^{j-1}=\phi^{j-1}\circ\phi_j -\phi^{j-1}$, it follows, for any $j\ge 2$,
\begin{align*}
\|\mathsf{W}_1(\phi^{j}-\phi^{j-1})\|_{r_j,s_j,\mathscr D_j}&\leby{VIPRel} \|\mathsf{W}_1D\phi^{j-1}\mathsf{W}_{j-1}^{-1}\|_{{r_{j-1}/2,s_{j-1}/2,}\mathscr D_{j-1}}\|\mathsf{W}_{j-1}\mathsf{W}_{j}^{-1}\| \|\mathsf{W}_j(\phi_j-\id)\|_{r_j,s_j,\mathscr D_j}\\
	&\overset{\equ{estfin2Ext03v501}+\equ{hihih}+\equ{DevphiJ}}{\le} 4\x^{2\n}\x^{(m+1)(j-1)}\;,
\end{align*}
so that
\beq{EstFrWhitn}
\sum_{j\ge 2}r_j^{-\wt m}\|\mathsf{W}_1(\phi^{j}-\phi^{j-1})\|_{r_j,s_j,\mathscr D_j}\le 4\x^{2\n}r_1^{-\wt m}\sum_{j\ge 2}\x^{(m+1)(j-1)}<\infty,
\eeq
from which we conclude that $\phi_*\in C^{\wt m}_W(\mathscr D_*\times\tn)$, and, in particular, 
\beq{Estvstaustar}
\sup_{\mathscr D_*\times \tn}\max\left\{|\mathsf{W}_1(\phi_*-\id)|,\;\|\dpr_x(u_*-\id)\|\right\}\le 8\x^{2\n}\;.
\eeq
%%$$
%%\sup_{\mathscr D_*\times \tn}|u_*-\id|\le \frac{\vae}{\a }\s_0^{l-\n},\qquad \sup_{\mathscr D_*\times \tn}\max\left\{\s^{-1}|u_*-\id|,\;\|\dpr_x u_*-\Id\|\right\}\le \frac{\mathsf{K}\vae}{\a^2}\s_0^{l-\n}
%%$$
Moreover,
$$
\sup_{y_*\in \mathscr D_*}\sum_{j\ge 2}s_j^{-\wt m \n}\sup_{\torus^d_{s_j}}|\mathsf{W}_1(\phi^{j}(y_*,\cdot)-\phi^{j-1}(y_*,\cdot))|\ltby{EstFrWhitn}\infty,
$$
which implies that $\phi_*$ is $C^{\wt m \n}$ in the angle variable \ie for any given $y_*\in\mathscr D_*$, the map $\phi_*(y_*,\cdot)\colon x\longmapsto \phi_*(y,x)$ is $C^{\wt m \n}(\tn)\subset C^{1}(\tn)$.

\noi
Now, by Lemma~\ref{BernsMos}, it follows that the sequence $\mathcal{H}_{j}$ converges in the $C^l$--topology uniformly to $H$ on $\rn\times\tn$. Thus, letting $j\rightarrow\infty$ in $K_j=\mathcal{H}_j\circ \phi^j-P_j$ yields
%converges uniformly to a function $K_*\in C^{\wh m}_W(\mathscr D_*)$ and, consequently, 
\equ{EqkStGSt} and \equ{EqHStfiStKst}.
\qed

\subsection{Completion of the Proof of Theorem~\ref{SecStatem}}
Choose  $\b\coloneqq \wt m$. Then, one checks easily that \equ{smcEAr0v2} implies \equ{SmaLConD} and, therefore, Lemmata~\ref{IteKAM} and \ref{ConvShc} hold. Thus, the map $G_0\coloneqq ({\dpr_y K_0}_{|B_{\tilde{r}_1/4}(\mathscr  D_{\a})})^{-1}\circ \dpr_y K$ is well-defined on $B_{\tilde{r}_0}(\mathscr D_{\a})$ and satisfies
\beq{G0DomEst}
G_0(B_{\tilde{r}_0}(\mathscr D_{\a}))\subset B_{\tilde{r}_1/2}(\mathscr D_{0})\;,\quad \max\left\{\|G_0-\id\|_{\tilde{r}_0,\mathscr D_\a},\;\x_0\|\dpr_y G_0-\uno_d\|_{\tilde{r}_0,\mathscr D_\a}\right\}\le 2C_1\eta\x_0^{l-1},
\eeq
 where $K_0\coloneqq \mathcal{K}_0$ and $\tilde{r}_0\coloneqq \tilde{r}_1/(16d\eta)$. Indeed, fix $y_0\in\mathscr D_\a$ and consider the auxiliary function $f\colon B_{\tilde{r}_1/4}(y_{0})\times B_{\tilde{r}_0}(y_0)\ni(y,z)\longmapsto \dpr_y K_0(y)-\dpr_y K(z)$. %, where $\hat{r}_0\coloneqq2\mathsf{T}\tilde{r}_0$. 
 Then, for any $(y,z)\in B_{\tilde{r}_1/4}(y_{0})\times B_{\tilde{r}_0}(y_{0})$
 \begin{align*}
  \|\uno_d-T(y_0)f_y(y,z)\|&\le \|T(y_0)\|\|\dpr_y^2(K_0-K)(y_0)+(\dpr_y^2 K_0(y)-\dpr_y^2 K_0(y_0))\|\\
  	&\leby{AnalAppro}\mathsf{T}(C_1\mathsf{K}\x_0^{l-2}+4dC_1\mathsf{K}\x_0^{-1}\tilde{r}_1/4){\le 2^{4l-7}C_1(\eta\s_0^{l-2}+\frac{\a}{\mathsf{K}}\s_0^\t)}\leby{smcEAr0v2} \su2,
 \end{align*}
 and
\begin{align*}
2\|T(y_0)\||f(y_0,z)|&\le 2\mathsf{T}|\dpr_y(K_0-K)(y_0)+(\dpr_y K(y_0)-\dpr_y K(z))|\\
	&\leby{AnalAppro} 2\mathsf{T}(C_1\mathsf{K}\x_0^{l-1}+d\mathsf{K}\tilde{r}_0)< 4d\eta\tilde{r}_0= \tilde{r}_1/4.
\end{align*}
Thus, by Lemma~\ref{IFTLem}, $G_0= (\dpr_y K_0)^{-1}\circ \dpr_y K$ is well-defined\footnote{In fact, the graph of $G_0$ is precisely the set of solutions of the equation $f(y,z)=0$.} on $B_{\tilde{r}_0}(\mathscr D_{\a})$ and the first part of \equ{G0DomEst} holds and we now prove its second part. In fact, for any  $y\in B_{\tilde{r}_0}(y_{0})$,
\begin{align*}
|G_0(y)-y|&=|(\dpr_{y}K_0)^{-1}( K_y(y))-(\dpr_{y}K_0)^{-1}\left(K_y(y)+\dpr_y (K_0-K)(y)\right)|\\
	   % &\le \dst\int_0^1\|\dpr_{y'}((\dpr_{y'}K')^{-1})(K_y(y)+t\vae\wt K_{y'}(y))\|dt\;\|\vae\wt K_{y'}\|_{\check{r},\mathscr D_\a}\\
	    &\le \|(\dpr_{y}^2K_0)^{-1}\|_{\tilde{r}_1/4,\mathscr D_\a}\|\dpr_y (K_0-K)\|_{C^l}\leby{AnalAppro} 2\mathsf{T}C_1\mathsf{K}\x_0^{l-1}
	  %  &=2^{2l-1}C_1\eta \s_0^{l-1}%\textcolor{red}{\ltby{smcEAr0v2} ???}
	    \;.
\end{align*}
Moreover,
$$
\|\dpr_y G_0-\uno_d\|_{\tilde{r}_0,\mathscr D_\a}=\|(\uno_d+T\dpr_y^2(K_0-K))^{-1}-\uno_d\|_{\tilde{r}_0,\mathscr D_\a}\le 2\mathsf{T}\|\dpr_y^2(K_0-K)\|_{\tilde{r}_0,\mathscr D_\a}\leby{AnalAppro} 2C_1\eta\x_0^{l-2},
$$
which completes the proof of \equ{G0DomEst}. 

\noi
Now, observe that $G^*= G_*\circ G_0 $ is well-defined by \equ{G0DomEst}. %where $G_0\coloneqq (\dpr_y \mathcal{K}_0)^{-1}\circ \dpr_y K$, 
Thus, $\mathscr D_0\bigcap B_{\tilde{r}_1/4}(\mathscr D_\a)=G_0(\mathscr D_\a)$ and, therefore, denoting $ G^*(\mathscr D_\a)$ again by ${ \mathscr D_*}$, the relations \equ{conjCaneq00v2} and \equ{conjCaneq0v2} then follows. 
%%
\begin{comment}
Next, we estimate $\phi_*$. We have, for any $i\ge 2$,
\begin{align*}
\|\mathsf{W}_1(\phi^i-\id)\|_{r_{i},s_{i},\mathscr{D}_{i}}&\le \|\mathsf{W}_1(\phi^{i-1}\circ\phi_i-\phi_i)\|_{r_{i},s_{i},\mathscr{D}_{i}}+\|\mathsf{W}_1(\phi_i-\id)\|_{r_{i},s_{i},\mathscr{D}_{i}}\\
&\le \|\mathsf{W}_1(\phi^{i-1}-\id)\|_{r_{i-1},s_{i-1},\mathscr{D}_{i-1}}+ \left(\dst\prod_{j=1}^{i-1}\|\mathsf{W}_{j}\mathsf{W}_{j+1}^{-1}\| \right) \|\mathsf{W}_{i}(\phi_i-\id)\|_{r_{i},s_{i},\mathscr{D}_{i}}\\
&\overset{\equ{estfin2Ext03v501}+\equ{DevphiJ}}{\le} \|\mathsf{W}_1(\phi^{i-1}-\id)\|_{r_{i-1},s_{i-1},\mathscr{D}_{i-1}}+ \frac{C_4\mathsf{K}\vae\s_0^{l-2\n}}{\a^2}\;\x^{(l-2\n+1)(i-1)}\;,
\end{align*}
 when iterated,  yields
\begin{align*}
\|\mathsf{W}_1(\phi^i-\id)\|_{r_i,s_i,\mathscr{D}_i}&\le \frac{C_4\mathsf{K}\vae\s_0^{l-2\n}}{\a^2}\dst\sum_{j\ge 1}\x^{(l-2\n+1)(i-1)}\le 2C_4\frac{\mathsf{K}\vae}{\a^2}\s_0^{l-2\n}\ltby{smcEAr0v2}\frac{\mathsf{K}\vae}{\a^2}\ltby{smcEAr0v2}1\,.\marginnote{$2C_4 \s_0^{l-2\n}\le 1$}
\end{align*}
\end{comment}
Observe that \equ{Estvstaustar} implies \equ{estArnTrExtv2}.

\noi
Next, we prove \equ{NormGstrThtv2}. Set $G^0\coloneqq G_0$, $G^{-1}\coloneqq \id$ and $\mathscr D_{-1}\coloneqq \mathscr D_\a$. Then, for any $j\ge 0$,
\begin{align*}
\|G^{j}-\id\|_{\mathscr D_{0}}&= \sum_{i\ge 0}\|G^{i+1}-G^i\|_{\mathscr D_{0}}\le \sum_{i\ge 0}\|G_{i+1}\circ G^i-G^i\|_{\mathscr D_{0}}\\
		&= \sum_{i\ge 0}\|G_{i+1}-\id\|_{\mathscr D_{i}}\le \sum_{i\ge 0}\|G_{i+1}-\id\|_{\tilde{r}_{i+1},\mathscr D_{i}}\nonumber\\
		&\overset{\equ{estGiidv201}+\equ{G0DomEst}}{\le} 2^{2l-1}C_1\eta\s_0^{l-1} +2\tilde{r}_1\x^{2\n}\leby{smcEAr0v2}\a\s_0^\t/\eta^2,
 %\label{DstrClToDj}
\end{align*}
then, letting $j\to\infty$ yields the first part of \equ{NormGstrThtv2}.

\noi
Next, we show that $\|G^*-\id\|_{L,\mathscr D_{\a}}<1$, which will imply that\footnote{See  \cite[Proposition II.2.]{zehnder2010lectures}.} $G^*\colon \mathscr D_{\a}\overset{onto}{\longrightarrow}\mathscr D_*$ is a lipeomorphism. Indeed, for any $j\ge 0$, we have
\begin{align*}
\| G^j-\id\|_{L,\mathscr D_{\a}}+1&= \| (G_j-\id)\circ G^{j-1}+(G^{j-1}-\id)\|_{L,\mathscr D_{\a}}+1\\
  &\le\| G_j-\id\|_{L,G^{j-1}(\mathscr D_{\a})}\| G^{j-1}\|_{L,\mathscr D_{\a}}+\| G^{j-1}-\id\|_{L,\mathscr D_{\a}}+1\\
  &\le\| G_j-\id\|_{L,G^{j-1}(\mathscr D_{\a})}(\|G^{j-1}-\id\|_{L,\mathscr D_{\a}}+1)+\| G^{j-1}-\id\|_{L,\mathscr D_{\a}}+1\\
 % &= (\| G_j-\id\|_{L,\mathscr D_{j-1}}+1)(\| G^{j-1}-\id\|_{L,\mathscr D_{\a}}+1)\\
  &\le (\|\dpr_z G_j-\uno_d\|_{\tilde{r}_j/2,\mathscr D_{j-1}}+1)(\| G^{j-1}-\id\|_{L,\mathscr D_{\a}}+1)\\
	%  &\le \|\dpr_z G_j-\uno_d\|_{\tilde{r}_j,\mathscr D_{j-1}}(\|\dpr_z G^{j-1}-\uno_d \|_{\check{r}_j,\mathscr D_{\a}}+1)+\|\dpr_z G^{j-1}-\uno_d\|_{\check{r}_j,\mathscr D_{\a}}\\
	%  &= (\|\dpr_z G_j-\uno_d\|_{\tilde{r}_j,\mathscr D_{j-1}}+1)(\|\dpr_z G^{j-1}-\uno_d \|_{\check{r}_j,\mathscr D_{\a}}+1)-1\\
 % &\overset{\equ{estGidevidv2}+\equ{estG1devidv2}}{\le} (64d\eta_0\s_{j-1}^{\n+d}|\vae|^{2^{j-1}}\mathsf{L}_{j-1}+1)(\| G^{j-1}-\id\|_{L,\mathscr D_{\a}}+1)
\end{align*}
which iterated and using Cauchy's estimate leads to\footnote{Recall that $\ex^x-1\le x\ex^x\;,\ \forall\; x\ge 0$.}
\begin{align*}
\| G^j-\uno_d\|_{L,\mathscr D_{\a}} &\le -1+\dst\prod_{i=0}^\infty (\|\dpr_z G_j-\uno_d\|_{\tilde{r}_i/2,\mathscr D_{i-1}}+1)\\
     &\overset{\equ{estGiidv201}+\equ{G0DomEst}}{\le} -1+\exp\left(2^{2l-3}C_1\eta\s_0^{l-2}+ \sum_{i=1}^\infty \x^{2\n}\x^{m(j-1)}\right)\\
     &\le -1+\exp\left(2^{2l-3}C_1\eta\s_0^{l-2}+2\x^{2\n}\right)\ltby{smcEAr0v2}1/2\;.
%\label{LipGUpjInf1}
\end{align*}
Thus, letting $j\to\infty$, we get that $G^*$ is Lipschitz continuous and \equ{NormGstrThtv2} is proven. Observe that \equ{Estvstaustar} implies \equ{estArnTrExtv2}. 
%
\begin{comment}
 Let us now prove the bound on $\dpr_x u_*-\uno_d$   in \equ{estArnTrExtv2}.
 For, set
$$\tilde{u}_j\coloneqq u-\id,\qquad U^j\coloneqq\dpr_x\pi_2\phi^j=(\uno_d+\dpr_x\tilde{u}_1)\cdots(\uno_d+\dpr_x\tilde{u}_j)
%,\qquad \tilde{u}_*\coloneqq \lim \tilde{u}_j
.
$$
 Then, for any $j\ge 1$, we have
\begin{align*}
 \|U^j\|_{s_{j}}&\le (1+\|\dpr_x\tilde{u}_1\|_{s_1})\cdots(1+\|\dpr_x\tilde{u}_j\|_{s_j})\\
 	&\stackrel{\equ{estfin2Ext03v501}}{\le}\exp\left(\dst\sum_{k\ge 1}C_4\r\s_0^l\;\x^{(l-2\n)(j-1)}\right)\le \exp\left(2 C_4\r\s_0^l\right)\ltby{smcEAr0v2}\ex,
\end{align*}
 so that
 $$
 \|U^{j+1}-U^j\|_{s_*}=\|U^{j}\cdot(\uno_d+\dpr_x\tilde{u}_{j+1})-U^j\|_{s_*}\le \|U^j\|_{s_{j+1}}\|\dpr_x\tilde{u}_{j+1}\|_{s_{j+1}} \stackrel{\equ{estfin2Ext03v501}}{\le}\ex C_4\r\s_0^l\x^{(l-2\n)(j-1)},
 $$
 which implies
 $$
  \|U^{j}-\uno_d\|_{s_*}\le 2\ex C_4\r\s_0^l\leby{smcEAr0v2}\frac{\mathsf{K}\vae}{\a^2}\s^{l-2\n}  \leby{smcEAr0v2}\su2
 $$
 and then letting $j\rightarrow\infty$, we get the estimate on $\dpr_x u_*-\uno_d$.
\end{comment}
%% 
 Next, observe that, thanks to \cite[Theorem~6.2.2, page~148]{koudjinan2019quantitative} (see also \cite{Chierrchia2019StructKolm}), \equ{NormGstrThtv2} yields \equ{MesChPin}.
 
 \noi
 Finally, we show that each KAM torus, as a graph, is of class $C^\n(\tn)$. Set\footnote{By \equ{estfin2Ext03v501}, for any given $y\in D_{2r_j}(\mathscr D_j)$, the map $x\longmapsto u_j(y,x)$ is a real--analytic diffeomorphism from $\torus^d_{s_j}$ onto its image, and we denote by $u_j^{-1}(y,x)$ the value of its inverse  at $(y,x)$.} $\f_j(y,x)\coloneqq v_j(y,u_j^{-1}(y,x))$, $\f^1\coloneqq \f_1$ and $\f^{j+1}(y,x)\coloneqq \f^j(\f_{j+1}(y,x),x)$,  $j\ge 1$. By definition, $\f_j(y,x)\eqby{ArnTraKamExtv5} y+\dpr_x g_j(y,x)$ so that, for any $j\ge 2$,
\begin{align}
 \max\left\{\|\f_j-\pi_1\|_{r_j,s_j,\mathscr D_{j}},\;r_{j-1}\|\dpr_y\f_j-\uno_d\|_{r_j,s_j,\mathscr D_{j}}\right\}&\stackrel{\equ{GrafCNu}+\equ{calPjxijmoin1}}{\le} (1+2\r)C_0\frac{\|\mathcal{P}^j\|_{r_{j-1},s_{j-1},\mathscr D_{j-1}}}{\a\s_{j-1}^\n}\nonumber\\
    &\stackrel{\equ{ConvDfiJ104g}}{\le} C_2\r\s_0^{l}\frac{\mathsf{K}}{\vae}\x^{(l-2\n)j-2l}r_{j-1}\nonumber\\
    &\eqqcolon \z_1 \x^{(l-2\n)j}\x^{\n(j-1)}. \label{fiInj0}
\end{align}
%where $\pi_1$ denotes the projection on the ``$y$''--component. 
Now, observe that
$
\dpr_y\f^j=\dpr_y\f_1\cdot \dpr_y\f_2\cdots \dpr_y\f_j.
$ 
Thus, by the usual telescoping argument, \equ{fiInj0} yields
$$
\z_2\coloneqq\sup_{j\ge 2}\|\dpr_y\f^j\|_{r_j,s_j,\mathscr D_{j}}<\infty.
$$
Therefore,
\begin{align*}
\sup_{y_*\in \mathscr D_*}\sum_{j\ge 2}s_{j+1}^{-\n}\sup_{\torus^d_{s_{j+1}}}|\f^{j+1}(y_*,\cdot)-\f^{j}(y_*,\cdot)|&\leby{estGiidv201} \sum_{j\ge 2}s_{j+1}^{-\n}\|\f^{j+1}-\f^{j}\|_{\tilde{r}_{j+1},s_{j+1},\mathscr D_{j+1}}\\
   &\le \z_2\sum_{j\ge 2}s_{j+1}^{-\n}\|\f_{j+1}-\pi_1\|_{r_{j+1},s_{j+1},\mathscr D_{j+1}}\\
    &\leby{fiInj0} \z_1\z_2s_1^{-\n}\sum_{j\ge 2}\x^{(l-2\n)j}<\infty.
\end{align*}
Hence, $\f^j$ converges uniformly on $\mathscr D_*\times\tn$ to the the map $\f_*$ defined by $\f_*(y_*,x)\coloneqq v_*(u_*^{-1}(y_*,x),x)$. %$=\lim_{j\to\infty}\f^j(y_*,x)$. 
In particular, $\f_*\in C^0(\mathscr D_*\times\tn)$ and, for any given $y_*\in\mathscr D_*$, the map $x\longmapsto\f_*(y_*,x)$ is $C^{\n}(\tn)$ and its graph is precisely the KAM torus $\phi_*(y_*,\tn)$.
 \qed
%\frac{\|\mathcal{P}^j\|_{r_j,s_j,\mathscr D_j}}{{\a}r_j}\r^{-1}&\le C_2\s_0^{l}\frac{\mathsf{K}}{\vae}\x^{(l-2\n)j-2l}\le C_2\s_0^{l}\frac{\mathsf{K}}{\vae}\x^{-2l} \leby{SmaLConD} 1\;,\label{ConvDfiJ104g}
\appendix
\section*{Appendix}
\addcontentsline{toc}{section}{Appendices}
\setcounter{section}{0}
\renewcommand{\thesection}{\Alph{section}}

\appA{Reminders}

\subsection{Classical estimates (Cauchy, Fourier, Cohomological Equation)} 
\lemtwo{Cau}{\rm \cite{CC95}} 
{\bf 1.} Let $p\in \natural,\,r,s>0, y_0\in \cn$ and $f$ a real--analytic function $D_{r,s}(y_0)$ with 
$$
\|f\|_{r,s}\coloneqq \dst\sup_{D_{r,s}(y_0)}|f|.
$$
%$D\subset \cn$ be an open domain and $f\colon D_\r(D)\to \cn$ be analytic with bounded sup--norm $\|f\|_\r$, for some $\r>0$. 
Then,\\
{\bf (i)} For any multi--index $(l,k)\in \natural^d\times\natural^d$ with $|l|_1+|k|_1\le p$ and for any $0<r'<r,\, 0<s'<s$,\footnote{As usual, $\dpr_y^l\coloneqq \frac{\dpr^{|l|_1}}{\dpr y_1^{l_1}\cdots\dpr y_d^{l_d}},\, \forall\, y\in\rn,\, l\in\zn $.\label{notDevPart}}
\[\|\partial_{y}^l \partial_{x}^k f\|_{r',s'}\leq p!\; \|f\|_{r,s}(r-r')^{|l|_1}(s-s')^{|k|_1}.\]%\frac{\|f\|_\r}{\s}.\]
%\elem
%In the next lemma, we recall some properties of the Fourier's coefficients of an analytic function.
%\lem{fce}
%Let $N\in \natural,\,f\in \mathcal{A}_{r,s,h,d},\, 0<\s < s$ with $N>\frac{d-1}{\s}$. Then
%\beqano
{\bf (ii)}  For any $ k\in \zn$ and any $y\in D_r(y_0)$ %and $\o\in \O_{\a,h}$,
$$|f_k(y)|\leq \ex^{-|k|_1 s}\|f\|_{r,s}.
$$
%%%%{\bf (iii)}  For any $ 0<\s < s$ and any $N>\frac{d-1}{\s},$ 
%%%%$$
%%%%\|f-T_N f\|_{r,s-\s}\leq 4^d C_2 N^d \ex^{-N\s}\|f\|_{r,s}.
%%%%$$
%\eeqano
{\bf 2.} Let $p\in \natural,\,\o\in \dst\D_\a^\t$ and $f\in \mathcal{A}_{r,s}$ %$f\colon \tn\to \real$ smooth with bounded holomorphic extension to $\torus^d_s$ 
and $\average{f}=0$. Then, for any $0<\s<s$, the system
\[D_\o g=f,\quad\average{g}=0\]
has a unique solution in $ \mathcal{A}_{r,s-\s}$  %and there exist constants $C_p=C_p(d,\t)\ge 1$ and $k_p=k_p(d,\t)\ge 1$ 
such that for any multi--index $k\in \natural^d$ with $|k|_1= l$
\[\|\partial_x^k g\|_{r,s-\s}\leq C_l\frac{\|f\|_{s}}{\a} \s^{-(\t+l)} ,\]
%In particular, one can take $\bar B_1=C_0$ 
%and $k_1=\t$ and this value of $k_1$ is the best known 
where $C_l\coloneqq 2^{d+1-(\t+l)}\sqrt{\G(2(\t+l)+1)}$
 (see \cite{RH75,CC95}).
\elem
\subsection{Implicit and Inverse function Theorems}
Firstly, we recall the classical implicit function Theorem, in a quantitative framework.
\lemtwo{IFTLem}{\cite{CL12}}
Let $ r,s>0,\, n,m\in \natural,\, (y_0,x_0)\in \complex^n\times\complex^m$ and\footnote{Here, $D^n_r(z_0)$ denotes the ball in $\complex^n$ centered at $z_0$ and with radius $r$.} %\footnote{Let us point out that any other norm (different!) may be used on $\complex^n,\, \complex^m$ and $\complex^{n+m}$.}
\[F\colon (y,x)\in D^n_r(y_0)\times D^m_s(x_0)\subset \complex^{n+m}\mapsto F(y,x)\in\complex^n\]
be continuous with continuous Jacobian matrix $F_y$. Assume that $F_y(y_0,x_0)$ is invertible with inverse $T\coloneqq F_y(y_0,x_0)^{-1}$ %and that there exist $c,c'>0$ with $c+c'=1$
 such that
\beq{HypIFT}
\sup_{D^n_r(y_0)\times D^m_s(x_0)}\|\uno_n-TF_y(y,x)\|\leq c<1 \quad \mbox{and}\quad \sup_{ D^m_s(x_0)}|F(y_0,\cdot)|\leq \frac{(1-c) r}{\|T\|}.
\eeq
Then, there exists a unique continuous function $ g\colon D^m_s(x_0)\to D^n_r(y_0)$ such that the following are equivalent
\begin{itemize}
\item[$(i)$] $(y,x)\in D^n_r(y_0)\times D^m_s(x_0)$ and $F(y,x)=0$;
\item[$(ii)$] $x\in D^m_s(x_0)$ and $y=g(x)$.
\end{itemize}
Moreover, $g$ satisfies
\beq{EstIFT}
\sup_{D^m_s(x_0)}|g-y_0|\leq \frac{\|T\|}{1-c}\sup_{D^m_s(x_0)}|F(y_0,\cdot)|.
\eeq
\elem

\appB{Outline of the proof of Lemma~\ref{lem:1Extv5Simpl01}\label{appC}}
Here, we aim to sketch the proof of the general KAM step. We refer the reader to \cite{Chierrchia2019StructKolm, koudjinan2019quantitative} for more details.

\noi
\Giu
{\bf Step 1: Construction of the Arnold's transformation} The symplectomorphism $\phi'$ is  generated by the real--analytic map $y'\cdot x+  g(y',x)$ \ie
\beq{ArnTraKamExtv5}
\phi'\colon \left\{\begin{aligned}
y  &=y'+  g_x(y',x)\\
x' &=x+  g_{y'}(y',x)\, ,
\end{aligned}
\right.
\eeq
in such a way that
\beq{ArnH1Extv5}
\left\{
\begin{aligned}
& H':= H\circ \phi'=K'+ P'\quad\mbox{on }D_{r_1,s_1}(\mathscr D_\sharp')\,,\\
& \det\dpr_{y'}^2 K'(y')\not=0\,,\quad\qquad\quad\;\forall\; y'\in \mathscr D_\sharp'\,,\\
& \dpr_{y'} K'(\mathscr D_\sharp')= \dpr_{y} K(\mathscr D_\sharp)\,,
\end{aligned}
\right.
\eeq
with 
\beq{ArnDefPsExtv5}
\left\{
\begin{aligned}
%K'    &\coloneqq K(y')+  P_0(y')\eqqcolon K(y')+  \wt K(y')\\
&P'(y',x')\coloneqq P_+(y',\f(y',x'))\;,\qquad P_+\coloneqq P^\ppu+P^\ppd+ P^\ppt\;,\qquad P^\ppt \coloneqq  P-\wh P\;,\\
&P^\ppu \coloneqq \dst\int^1_0(1-t)K_{yy}(  t g_x)\cdot g_x\cdot g_x dt\;,\qquad  P^\ppd \coloneqq \dst\int_0^1P_y(y'+  t g_x,x)\cdot g_x dt\; ,
\end{aligned}
\right.
\eeq
where $\f(y',\cdot)$ is the inverse of the map $x' \mapsto x+  g_{y'}(y',x)$ and $\wh P$ is the approximation of $P$ given by \cite[Theorem~7.2]{russmann2001invariant}.\footnote{With the choices $\b_1=\cdots=\beta_d=1/2$, $T=\k$ and $\d=\r\leby{DefNArnExt1v501} 1/4$.} Moreover,
\beq{CondHomEqArn1v5}
K_y(y')\cdot n\not= 0, \quad \forall\; 0<|n|_1\leq \k,\quad \forall\; y'\in D_{r_1}(\mathscr D_\sharp')\quad  \left(\subset D_{r}(\mathscr D_\sharp)\right),
\eeq
and the generating function a $g$ is given by
 \beq{HomEqArn1v5}
 g(y',x)\coloneqq \dst\sum_{0<|n|_1\leq \k} \frac{-\wh P_n(y')}{iK_y(y')\cdot n}\ex^{in\cdot x}.
 \eeq
%\subsubsection{quant estimates}
{\bf Step 2} Now, we provide the construction performed in {\bf Step 1} with quantitative estimate. 
First of all, notice that\footnote{By definition of $\wh P$, see \cite[Theorem~7.2]{russmann2001invariant}.}
\beq{RussemanTroncation}
\|P-\wh P\|_{r,\bar{s},\mathscr D_\sharp}\le 2\r \vae\,, \qquad \|\wh P\|_{r,\bar{s}}\le \|P\|_{r,s,\mathscr D_\sharp}+\|P-\wh P\|_{r,\bar{s},\mathscr D_\sharp}\le(1+2\r)\vae.
\eeq
Observe also that for any $\mathsf{y}\in\mathscr D_\sharp$, $0<|n|_1\leq \k$ and $y'\in D_{\bar{r}}(\mathsf{y})$,
\begin{align}
|K_y(y')\cdot n|&\ge \frac{\a}{2|n|_1^\t}.\label{ArnExtDiopCondExtv5}
\end{align}
Now, using Lemma~\ref{Cau}--$2.$ and \equ{RussemanTroncation}, we get %(recall that $\bar{r}\leq \frac{r}{2}$)
\beq{GrafCNu}
\begin{aligned}
&\|g\|_{\bar{r},\bar{s},\mathscr D_\sharp} \le \mathsf{C}_0 \frac{(1+2\r)\vae}{\a} \s^{-\t}\,, \quad \|g_x\|_{\bar{r},\bar{s},\mathscr D_\sharp} \le \mathsf{C}_0 \frac{(1+2\r)\vae}{\a} \s^{-(\t+1)}\,, \\
& \max\{\|\dpr_{y'}g\|_{\bar{r},\bar{s},\mathscr D_\sharp},\ \s\|\dpr^2_{y'x}g\|_{\bar{r},\bar{s},\mathscr D_\sharp}\;,\ \s^2\|\dpr^3_{y'xx}g\|_{\bar{r},\bar{s},\mathscr D_\sharp}\}\le \ovl{\mathsf{L}},
\end{aligned}
\eeq
where
$$
\ovl{\mathsf{L}}\coloneqq 2\mathsf{C}_0 \frac{(1+2\r)\vae}{\a r} \s^{-\t}\;.
$$
We have
\begin{align*}
&\|\dpr_{y'}\wt K\|_{r/2,\mathscr D_\sharp}\le  \frac{2\vae}{r} \;,\qquad \|\dpr_{y'}^2\wt K\|_{r/2,\mathscr D_\sharp}\leq \frac{4\vae}{r^2}\leby{cond1ExtExtv501} \mathsf{K}\frac{\bar{\s}}{3}\;.
\end{align*}
\noi
Next, we construct $\mathscr D_\sharp'$ %and uniqueness of $y_1$
 in \eqref{ArnH1Extv5}. 
 For, fix  $\mathsf{y}\in\mathscr D_\sharp$ %, $z\in B^d_{\bar{r}}(y_0)$
  and consider
\begin{align*}
F\colon D_{\check{r}}(\mathsf{y})\times D_{\tilde{r}}(\mathsf{y}) &\longrightarrow \qquad \cn\\
		(y,z)\quad &\longmapsto K_y(y)+  \wt K_{y'}(y)-K_y(z).
\end{align*}
Then, one checks easily that Lemma~\ref{IFTLem} applies. Thus, we get that $F^{-1}(\{0\})$ is given by the graph of a real--analytic map $G^{\mathsf{y}}\colon D_{\tilde{r}}(\mathsf{y})\to D_{\check{r}}(\mathsf{y})$. Afterwards, one checks that the pieces of the family $\{G^{\mathsf{y}}\}_{\mathsf{y}\in\mathscr D_\sharp}$ matches, yielding therefore a global map $G$ on $D_{\tilde{r}}(\mathscr D_\sharp)$ and that, in fact, $G$ is bi--real--analytic.\footnote{\ie an invertible real--analytic map whose inverse is real--analytic as well.} Next, one shows that the expression $(K_y+ \wt K_{y'})^{-1}\circ K_y$ defines a map on $D_{\tilde{r}}(\mathsf{y})$ by means of the Inversion Function Lemma~\ref{IFTLem}. As a consequence, we  get an explicit formula for $G$:
\beq{explForG}
G=(K_y+ \wt K_{y'})^{-1}\circ K_y\qquad \mbox{on}\quad D_{\tilde{r}}(\mathsf{y})\;,
\eeq
and $\mathscr D_\sharp'=G(\mathscr D_\sharp)$. 
The reminder of the proof then goes exactly as in \cite{Chierrchia2019StructKolm} (see also \cite{koudjinan2019quantitative}). %Lemma~\ref{lem:1}.}
\qed

\appC{Isotropicity Lemma\label{appC}}
\lemtwo{LagTor}{\cite{herman1988existence,herman1989inegalites,broer2009quasi}}%{{\cite{salamon2004kolmogorov}}}
Let $H\colon \mathcal{M}\coloneqq \rn\times\tn\to \real$ be a Hamiltonian  of class $C^1$ and $\phi: \tn\ni x\longmapsto \phi(x)\in \mathcal{M}$, a $C^1$--mapping. Assume that 
$$
\phi_H^t\circ \phi(x)= \phi(x+t\o), \qquad \forall\; x\in \tn,
$$
for some rationally independent\footnote{\ie for any $k\in \zn\setminus\{0\}$, $\o\cdot k\not=0$.} $\o\in\rn$. Then, the torus $\phi(\tn)$ is isotropic \ie $\mathit{i}^*\varpi\equiv 0$, where $\mathit{i}\colon \phi(\tn)\hookrightarrow\mathcal{M}$ is the inclusion map.%a Lagrangian submanifold of $\mathcal{M}$.
\elem

\noi
{\bf Acknowledgments.} %We are grateful to 
The author would like to thank Prof. Luigi Chierchia, U. Bessi and A. Bounemoura for  fruitful discussions and criticisms on an earlier draft of this manuscript. The author would like to express his gratitude to Prof. Luigi Chierchia and S. Luzzatto for their support. % useful discussions and valuable comments on the first draft of this paper. 
%Theauthor would like to express their gratitude to ??? for stimulatingdiscussions that helped improve the manuscript.
The work of this paper has been completed during the author visit at the Abdus Salam International Centre for Theoretical Physics in Trieste, Italy. 
%%%%%%%%%%%%
%\nocite{koudjinan2019quantitative}
\bibliography{BibtexDatabase}
\end{document}